\documentclass[11pt,twoside,letterpaper]{article}
\usepackage{amssymb}
\usepackage{amsfonts}
\usepackage{amsmath}
\usepackage{times,fancyhdr}
\usepackage[dvips]{graphicx}

\makeatletter
\def\cleardoublepage{\clearpage\if@twoside \ifodd\c@page\else    \hbox{}    \thispagestyle{empty}    \newpage    \if@twocolumn\hbox{}\newpage\fi\fi\fi}
\makeatother

\def\figurename{Figure}
\makeatletter
\renewcommand{\fnum@figure}[1]{\figurename~\thefigure.}
\makeatother
\def\tablename{Table}
\makeatletter
\renewcommand{\fnum@table}[1]{\tablename~\thetable.}
\makeatother
\input{tcilatex}
\begin{document}

\title{ {\vskip0.45in \bfseries\scshape}\textbf{A General Maximum Principle
for Mean-field Stochastic Differential Equations with Jump Processes}}
\author{\bfseries\itshape Mokhtar Hafayed\thanks{\emph{Laboratory of Applied
Mathematics, Biskra University, Po Box 145, Biskra 07000, Algeria. E-mail
address: hafayedmokhtar@yahoo.com}} , Syed Abbas\thanks{\emph{School of
Basic Sciences, Indian Institute of Technology Mandi, Mandi H.P. 175001
India. E-mail address: sabbas.iitk@gmail.com}}}
\maketitle

\noindent \textbf{Abstract.} \emph{In this paper, we investigate the optimal
control problems for stochastic differential equations (SDEs in short) of
mean-field type with jump processes. The control variable is allowed to
enter into both diffusion and jump terms. This stochastic maximum principle
differs from the classical one in the sense that here the first-order
adjoint equation turns out to be a linear mean-field backward SDE with
jumps, while the second-order adjoint equation remains the same as in Tang
and Li's stochastic maximum principle }\cite{tang}\emph{. Finally, for the
reader's convenience we give some analysis results used in this paper in the
Appendix.}$\smallskip \smallskip $

\noindent \textbf{Keywords:\ }\textit{Mean-field SDEs. Jump processes.
Optimal stochastic control. Maximum principle. Spike varition method.
McKean-Vlasov equations.}$\smallskip $


\section{\textbf{Introduction}}

\noindent In this paper we study stochastic optimal control for a system
governed by nonlinear SDEs of mean-field type, which is also called
McKean-Valasov equations, with jump processes:

\begin{equation}
\left\{ 
\begin{array}{c}
dx^{u}(t)=f\left( t,x^{u}(t),\mathbb{E(}x^{u}(t)),u(t)\right) dt+\sigma
\left( t,x^{u}(t),\mathbb{E(}x^{u}(t)),u(t)\right) dW(t)\smallskip
\smallskip  \\ 
+\dint_{\Theta }g\left( t,x^{u}(t_{-}),u(t),\theta \right) N\left( d\theta
,dt\right) ,\smallskip \smallskip  \\ 
\multicolumn{1}{l}{x^{u}(s)=\zeta ,}%
\end{array}%
\right.   \tag{1.1}  \label{1.1}
\end{equation}%
where the coefficients $f$ and $\sigma $ depend on the state of the solution
process as well as of its expected value and the initial time $s$ and the
initial state $\zeta $ of the system are fixed, $(W(t))_{t\in \left[ s,T%
\right] }$ is a standard $one-$dimentional Brownian motion and $N\left(
d\theta ,dt\right) $ is a Poisson martingale measure with characteristic $%
\mu \left( d\theta \right) dt.$ This mean-field jump diffusion processes are
obtained as the mean-square limit, when $n\rightarrow +\infty $ of a system
of interacting particles of the form%
\begin{eqnarray*}
dx_{n}^{j,u}(t) &=&f\left( t,x_{n}^{j,u}(t),\frac{1}{n}\sum%
\limits_{i=1}^{n}x_{n}^{i,u}(t),u(t)\right) dt \\
&&+\sigma \left( t,x_{n}^{j,u}(t),\frac{1}{n}\sum%
\limits_{i=1}^{n}x_{n}^{i,u}(t),u(t)\right) dW^{j}(t) \\
&&+\int_{\Theta }g\left( t,x_{n}^{j,u}(t_{-}),u(t),\theta \right)
N^{j}\left( d\theta ,dt\right) .
\end{eqnarray*}%
\medskip 

\noindent Our control problem consists in minimizing a cost functional of
the form:

\begin{equation}
J^{^{s,\zeta }}\left( u(\cdot )\right) =\mathbb{E}\left[ h(x^{u}(T),\mathbb{E%
}\left( x^{u}(T)\right) )+\int_{s}^{T}\mathfrak{\ell }(t,x^{u}(t),\mathbb{E}%
\left( x^{u}(t)\right) ,u(t))dt\right] .  \tag{1.2}  \label{1.2}
\end{equation}%
This cost functional is also of mean-field type, as the functions $h$ and $%
\ell $ depend on the marginal law of the state process through its expected
value.

\noindent An admissible control $u(\cdot )$\ is an $\mathcal{F}_{t}$-adapted
and square-integrable process with values in a nonempty subset $\mathbb{A}$
of $\mathbb{R}.$ We denote the set of all admissible controls by $\mathcal{U}
$. Any admissible control $u(\cdot )\in \mathcal{U}$ satisfying%
\begin{equation}
J^{^{s,\zeta }}\left( u^{\ast }(\cdot )\right) =\min_{u(\cdot )\in \mathcal{U%
}}J^{^{s,\zeta }}\left( u(\cdot )\right) ,  \tag{1.3}  \label{1.3}
\end{equation}%
is called an optimal control. The corresponding state process, solution of
SDE-(\ref{1.1}), is denoted by $x^{\ast }(\cdot )=x^{u^{\ast }}(\cdot ).$

The modern optimal control theory has been well developed since early 1960s,
when Pontryagin et al., \cite{pontry} published their work on the maximum
principle and Bellman \cite{belman} put forward the dynamic programming
method. The pioneering works on the stochastic maximum principle was written
by Kushner (\cite{kushner},\cite{kushner1}). Since then there have been a
lot of works on this subject, see for instance (\cite{peng},\cite{zhoo},\cite%
{bensoussan},\cite{cad},\cite{haussmann},\cite{hafayed1},\cite{wang}). Peng 
\cite{peng} obtained the optimality stochastic maximum principle for the
general case. A good account and an extensive list of references on
stochastic optimal control can be founded in Yong et al., \cite{yong}.$%
\smallskip $

\noindent The stochastic optimal control problems for jump processes has
been investigated by many authors, see for instance, (\cite{cad},\cite%
{framstad},\cite{oksindal1},\cite{rishel},\cite{tang},\cite{hafayed2},\cite%
{shi1},\cite{shi2},\cite{situ}). The stochastic maximum principle for jump
diffusion in general case, where The control domain need not be convex. and
the diffusion coefficient depends explicitly on the control variable, was
derived via spike variation method by Tang et al., \cite{tang}, extending
the Peng's stochastic maximum principle of optimality \cite{peng}. These
conditions are described in terms of two adjoint processes, which are linear
classical backward SDEs. The sufficient conditions for optimality was
obtained by Framstad et al., \cite{framstad}.$\smallskip $

\noindent Historically, the SDE of Mean-field type was introduced by Kac 
\cite{kac} in 1956 as a stochastic model for the Vlasov-kinetic equation of
plasma and the study of which was initiated by McKean \cite{mckean} in 1966.
Since then, many authors made contributions on SDEs of mean-field type and
applications, see for instance, (\cite{ahmed},\cite{buckdahn2},\cite%
{buckdahn1},\cite{borkar},\cite{kac},\cite{li},\cite{mayer},\cite{vere},\cite%
{shi5},\cite{yong1}). Mean- field stochastic maximum principle of
optimality\ was considered by many authors, see for instance (\cite%
{buckdahn1},\cite{li},\cite{mayer},\cite{yong1}). In Buckdahn et al., \cite%
{buckdahn2} the authors obtained mean-field backward stochastic differential
equations. In a recent paper by Buckdahn et al., \cite{buckdahn1}, the
maximum principle was introduced for a class of stochastic control problems
involving SDEs of mean-field type, where the authors obtained a stochastic
maximum principle differs from the classical one in the sense that the
first-order adjoint equation turns out to be a linear mean-field backward
SDE, while the second-order adjoint equation remains the same as in Peng's
stochastic maximum principle \cite{peng}. In Mayer-Brandis et al., \cite%
{mayer} a stochastic maximum principle of optimality for systems governed by
controlled It\^{o}-Levy process of mean-field type is proved by using
Malliavin calculus. The local maximum principle of optimality for Mean-field
stochastic control problem has been derived by Li \cite{li}. The
linear-quadratic optimal control problem for mean-field SDEs has been
studied by Yong \cite{yong1}.

Our purpose in this paper is to establish necessary conditions of optimality
for Mean-field SDEs with jumps processes, in which the coefficients of
diffusion depend on the state of the solution process as well as of its
expected value. Moreover, the cost functional is also of Mean-field type.
The proof of our main result is based on spike variation method. This
results is an extension of \textit{Theorem 2.1} in Buckdahn et al., \cite%
{buckdahn1} to the controlled mean-field SDEs with jump processes. To
streamline the presentation, we only consider the one dimensional case.

The rest of the paper is organized as follows. Section 2 begins with a
general formulation of a mean-field control problem with jump processes and
give the notations and assumptions used throughout the paper. In Sections 3
we prove our main result.

\section{\textbf{Assumptions and statement of the control problem}}

\noindent Let $(\Omega ,\mathcal{F},\left( \mathcal{F}_{t}\right) _{t\in %
\left[ 0,T\right] },\mathbb{P})$ be a fixed filtered probability space
equipped with a $\mathbb{P}-$completed right continuous filtration on which
a $d-$dimensional Brownian motion $W=\left( W(t)\right) _{t\in \left[ 0,T%
\right] }$ is defined. Let $\eta $ be a homogeneous $\left( \mathcal{F}%
_{t}\right) $-Poisson point process independent of $W$. We denote by $%
\widetilde{N}(d\theta ,dt)$ the random counting measure induced by $\eta $,
defined on $\Theta \times \mathbb{R}_{+}$, where $\Theta $ is a fixed
nonempty subset of $\mathbb{R}$ with its Borel $\sigma $-field $\mathcal{B}%
\left( \Theta \right) $. Further, let $\mu \left( d\theta \right) $ be the
local characteristic measure of $\eta $, i.e. $\mu \left( d\theta \right) $
is a $\sigma $-finite measure on $\left( \Theta ,\mathcal{B}\left( \Theta
\right) \right) $ with $\mu \left( \Theta \right) <+\infty $. We then define 
\begin{equation*}
N(d\theta ,dt)=\widetilde{N}(d\theta ,dt)-\mu \left( d\theta \right) dt,
\end{equation*}%
where $N$ is Poisson martingale measure on $\mathcal{B}\left( \Theta \right)
\times \mathcal{B}\left( \mathbb{R}_{+}\right) $ with local characteristics $%
\mu \left( d\theta \right) dt.$ We assume that $\left( \mathcal{F}%
_{t}\right) _{t\in \left[ 0,T\right] }$ is $\mathbb{P}-$augmentation of the
natural filtration $(\mathcal{F}_{t}^{(W,N)})_{t\in \left[ s,T\right] }$
defined as follows%
\begin{equation*}
\mathcal{F}_{t}^{(W,N)}=\sigma \left( W(r):s\leq r\leq t\right) \vee \sigma
\left( \int_{s}^{r}\int_{B}N(d\theta ,d\tau ):s\leq \tau \leq t,\text{ }B\in 
\mathcal{B}\left( \Theta \right) \right) \vee \mathcal{G},
\end{equation*}%
where $\mathcal{G}$ denotes the totality of $\mathbb{P}-$null sets, and $%
\sigma _{1}\vee \sigma _{2}$ denotes the $\sigma $-field generated by $%
\sigma _{1}\cup \sigma _{2}.$ \newline

\noindent \textbf{Basic notations. }For convenience, we will use the
following notations throughout the paper. Let $u(\cdot )\in \mathcal{U}$ be
an admissible control. For $\Phi =f,\sigma ,\ell :$

\begin{enumerate}
\item $\delta \Phi (t)=\Phi (t,x^{\ast }(t),\mathbb{E}(x^{\ast
}(t)),u(t))-\Phi (t,x^{\ast }(t),\mathbb{E}(x^{\ast }(t)),u^{\ast }(t)).$

\item $\Phi _{x}(t)=\dfrac{\partial \Phi }{\partial x}(t,x^{\ast }(t),%
\mathbb{E}(x^{\ast }(t)),u^{\ast }(t))$, $\Phi _{y}(t)=\dfrac{\partial \Phi 
}{\partial y}(t,x^{\ast }(t),\mathbb{E}(x^{\ast }(t)),u^{\ast }(t)).$

\item $g_{x}\left( t,\theta \right) =g_{x}\left( t,x(t_{-}),u(t),\theta
\right) ,$ $g_{xx}\left( t,\theta \right) =g_{xx}\left(
t,x(t_{-}),u(t),\theta \right) .$

\item $\Phi _{xx}(t)=\frac{\partial ^{2}\Phi }{\partial x^{2}}(t,x^{\ast
}(t),\mathbb{E}(x^{\ast }(t)),u^{\ast }(t)),$ $\Phi _{yy}(t)=\frac{\partial
^{2}\Phi }{\partial y^{2}}(t,x^{\ast }(t),\mathbb{E}(x^{\ast }(t)),u^{\ast
}(t)).$

\item $\Phi _{xy}(t)=\frac{\partial ^{2}\Phi }{\partial x\partial y}%
(t,x^{\ast }(t),\mathbb{E}(x^{\ast }(t)),u^{\ast }(t)).$

\item $\mathcal{L}_{t}(\Phi ,y)=\frac{1}{2}\Phi _{xx}(t,x^{\ast }(t),\mathbb{%
E}(x^{\ast }(t)),u^{\ast }(t))y^{2}$, $\mathcal{L}_{t,\theta }(g,y)=\frac{1}{%
2}g_{xx}(t,x^{\ast }(t),u^{\ast }(t),\theta )y^{2}.$

\item We denote by $\mathbf{I}_{A}$ the indicator function of $A$ and by $%
sgn(\cdot )$ the sign function.

\item We denote by $\mathbb{L}_{\mathcal{F}}^{2}\left( \left[ s,T\right] ;%
\mathbb{R}\right) =\left\{ \phi (\cdot ):=\phi (t,w)\text{ is an }\mathcal{F}%
_{t}-\text{adapted }\mathbb{R}-\text{valued}\right. $ $\left. \text{%
measurable process on }\left[ s,T\right] \text{ such that }\mathbb{E}\left(
\int_{s}^{T}\left\vert \phi (t)\right\vert ^{2}dt\right) <\infty \right\} ,$
and by

$\mathbb{M}_{\mathcal{F}}^{2}\left( \left[ s,T\right] ;\mathbb{R}\right)
=\left\{ \phi (\cdot ):=\phi (t,\theta ,w)\text{ is an }\mathcal{F}_{t}-%
\text{adapted }\mathbb{R}-\text{valued measurable}\right. $

$\left. \text{ process on }\left[ s,T\right] \times \Theta \text{ such that }%
\mathbb{E}\left( \int_{s}^{T}\int_{\Theta }\left\vert \phi (t,\theta
)\right\vert ^{2}\mu \left( d\theta \right) dt\right) <\infty \right\} .$

\item In what follows, $C$ and $\rho (\varepsilon )$ represents a generic
constants, which can be different from line to line.
\end{enumerate}

\medskip

\noindent \textbf{Basic assumptions.} Throughout this paper we assume the
following.

\begin{description}
\item[\textbf{(H1)}] The functions $f(t,x,y,u):\left[ s,T\right] \times 
\mathbb{R}\times \mathbb{R\times \mathbb{A}\rightarrow R},$ $\sigma
(t,x,y,u):\left[ s,T\right] \times \mathbb{R}\times \mathbb{R\times \mathbb{A%
}\rightarrow R}$, $\ell (t,x,y,u):\left[ s,T\right] \times \mathbb{R}\times 
\mathbb{R\times \mathbb{A}}\rightarrow \mathbb{R}$ and $h(x,y):\mathbb{R}%
\times \mathbb{R\rightarrow R}$ are twice continuously differentiable with
respect to $\left( x,y\right) $. Moreover, $f,\sigma ,h$ and $\ell $ and all
their derivatives up to second-order with respect to $\left( x,y\right) $
are continuous in $\left( x,y,u\right) $ and bounded.

\item[\textbf{(H2)}] The function $g:\left[ s,T\right] \times \mathbb{%
R\times }\mathbb{A\times }\Theta \rightarrow \mathbb{R}$ is twice
continuously differentiable in $x$, Moreover $g_{x}$ is continuous, $%
\sup_{\theta \in \Theta }\left\vert g_{x}(t,\theta )\right\vert <+\infty $
and there exists a constant $C>0$ such that%
\begin{equation}
\begin{array}{c}
\sup_{\theta \in \Theta }\left\vert g\left( t,x,u,\theta \right) -g\left(
t,x^{\prime },u,\theta \right) \right\vert +\sup_{\theta \in \Theta
}\left\vert g_{x}\left( t,x,u,\theta \right) -g_{x}\left( t,x^{\prime
},u,\theta \right) \right\vert \\ 
\leq C\left\vert x-x^{\prime }\right\vert%
\end{array}
\tag{2.1}  \label{2.1}
\end{equation}%
\begin{equation}
\sup_{\theta \in \Theta }\left\vert g\left( t,x,u,\theta \right) \right\vert
\leq C\left( 1+\left\vert x\right\vert \right) .  \tag{2.2}  \label{2.2}
\end{equation}
\end{description}

\noindent Under the above assumptions, the SDE-(\ref{1.1}) has a unique
strong solution $x^{u}(t)$ which is given by%
\begin{eqnarray*}
x^{u}(t) &=&\zeta +\int_{s}^{t}f\left( r,x^{u}(r),\mathbb{E}%
(x^{u}(r)),u(r)\right) dr+\int_{s}^{t}\sigma \left( r,x^{u}(r),\mathbb{E}%
(x^{u}(r)),u(r)\right) dW(r) \\
&&+\int_{s}^{t}\int_{\Theta }g\left( t,x^{u}(r_{-}),u(r),\theta \right)
N\left( d\theta ,dr\right) ,
\end{eqnarray*}%
and by standard arguments it is easy to show that for any $q>0$, it holds
that

\begin{equation}
\mathbb{E(}\sup_{t\in \left[ s,T\right] }\left\vert x^{u}(t)\right\vert
^{q})<C_{q},  \tag{2.3}  \label{2.3}
\end{equation}%
where $C_{q}$ is a constant depending only on $q$ and the functional $%
J^{s,\zeta }$ is well defined.

\noindent \textbf{Usual Hamiltonian. }We define the usual Hamiltonian
associated with the mean-field stochastic control problem (\ref{1.1})-(\ref%
{1.2}) as follows%
\begin{equation}
\begin{array}{l}
H\left( t,X,\mathbb{E}\left( X\right) ,u,\Psi (t),K(t),\gamma _{t}(\theta
)\right) =\Psi (t)f\left( t,X,\mathbb{E}\left( X\right) ,u\right) \smallskip
\smallskip \\ 
\text{ \ \ \ \ \ \ \ }+K(t)\sigma \left( t,X,\mathbb{E}\left( X\right)
,u\right) +\dint_{\Theta }\gamma _{t}(\theta )g\left( t,x(t),u(t),\theta
\right) \mu (d\theta )\smallskip \smallskip \\ 
\text{ \ \ \ \ \ \ \ }-\ell \left( t,X,\mathbb{E}\left( X\right) ,u\right) ,%
\end{array}
\tag{2.4}  \label{2.4}
\end{equation}%
where $(t,X,u)\in \lbrack s,T]\times \mathbb{R}\times \mathbb{A}$, $X$\ is a
random variable such that $X\in \mathbb{L}^{1}\left( \Omega ,\mathcal{F},%
\mathbb{R}\right) $ and $\left( \Psi (t),K(t),\gamma _{t}(\theta )\right)
\in \mathbb{R\times R\times R}$ given by equation \textit{(\ref{2.5}).}$%
\smallskip $

\noindent \textbf{Adjoint equations for mean-field SDEs with jump processes. 
}We introduce the adjoint equations involved in the stochastic maximum
principle for our control problem. The first-order adjoint equation turns
out to be a linear mean-field backward SDE with jump terms, while the
second-order adjoint equation remains the same as in Tang et al., \cite{tang}%
.

\noindent For any $u(\cdot )\in \mathcal{U}$ and the corresponding state
trajectory $x(\cdot )$, we define the first-order adjoint process $(\Psi
(\cdot ),K(\cdot ),\mathbf{\gamma }(\cdot ))$ and the second-order adjoint
process $(Q(\cdot ),R(\cdot ),\Gamma (\cdot ))$ as the ones satisfying the
following equations:

\begin{enumerate}
\item[\textbf{1.}] \textit{First-order adjoint equation: linear backward SDE
of mean-field type with jump processes}%
\begin{equation}
\left\{ 
\begin{array}{l}
d\Psi (t)=-\left\{ f_{x}\left( t\right) \Psi (t)+\mathbb{E}\left(
f_{y}^{\top }(t)\Psi (t)\right) \right. +\sigma _{x}\left( t\right)
K(t)\smallskip \smallskip \\ 
\text{ \ \ \ \ \ \ \ \ }+\mathbb{E}\left( \sigma _{y}(t)K(t)\right) +\ell
_{x}\left( t\right) +\mathbb{E}\left( \ell _{y}(t\right) )\smallskip
\smallskip \\ 
\text{ \ \ \ \ \ \ \ \ }+\left. \dint_{\Theta }g_{x}\left( t,\theta \right) 
\mathbf{\gamma }_{t}(\theta )\mu (d\theta )\right\} dt\smallskip \smallskip
\\ 
\text{ \ \ \ \ \ \ \ \ }+K(t)dW(t)+\dint_{\Theta }\mathbf{\gamma }%
_{t}(\theta )N(dt,d\theta )\smallskip \smallskip \\ 
\Psi (T)=-\left( h_{x}\left( x(T),\mathbb{E}(x(T)\right) +\mathbb{E}\left(
h_{y}\left( x(T),\mathbb{E}(x(T)\right) \right) \right) .%
\end{array}%
\right.  \tag{2.5}  \label{2.5}
\end{equation}

\item[\textbf{2.}] \textit{Second-order adjoint equation: classical linear
backward SDE with jump processes (see Tang et al., \cite{tang} equation
(2.23))}%
\begin{equation}
\left\{ 
\begin{array}{l}
dQ(t)=-\left\{ 2f_{x}\left( t\right) Q(t)+\sigma _{x}^{2}\left( t\right)
Q(t)+2\sigma _{x}\left( t\right) R(t)\right. \smallskip \smallskip \\ 
\text{ \ \ \ }+\dint_{\Theta }\left( \Gamma _{t}(\theta )+Q(t)\right) \left(
g_{x}\left( t,\theta \right) \right) ^{2}\mu (d\theta )+2\dint_{\Theta
}\Gamma _{t}(\theta )g_{x}\left( t,\theta \right) \mu (d\theta )\smallskip
\smallskip \\ 
\text{ \ \ \ }+\left. H_{xx}(t))\right\} dt+R(t)dW(t)+\dint_{\Theta }\Gamma
_{t}(\theta )N(d\theta ,dt)\smallskip \smallskip \\ 
Q(T)=-h_{xx}\left( x(T),\mathbb{E}(x(T))\right) .%
\end{array}%
\right.  \tag{2.6}  \label{2.6}
\end{equation}
\end{enumerate}

\noindent \textbf{Remark 2.1.}\ As it is well known that under conditions
(H1) and (H2) the first-order adjoint equation\ (\ref{2.5}) admits one and
only one $\mathcal{F}_{t}-$adapted solution pair $\left( \Psi (\cdot
),K(\cdot ),\mathbf{\gamma }(\cdot )\right) \in $\ $\mathbb{L}_{\mathcal{F}%
}^{2}\left( \left[ s,T\right] ;\mathbb{R}\right) $\ $\times \mathbb{L}_{%
\mathcal{F}}^{2}\left( \left[ s,T\right] ;\mathbb{R}\right) $\ $\times 
\mathbb{M}_{\mathcal{F}}^{2}\left( \left[ s,T\right] ;\mathbb{R}\right) $.
This equation reduces to the standard one as in (Tang et al., \cite{tang}\
equation (2.22)), when the coefficients not explicitly depend on the
expected value (or the marginal law) of the underlying diffusion process.
Also the second-order adjoint equation (\ref{2.6}) admits one and only one $%
\mathcal{F}_{t}-$adapted solution pair $\left( Q(\cdot ),R(\cdot ),\Gamma
(\cdot ,\theta )\right) \in \mathbb{L}_{\mathcal{F}}^{2}\left( \left[ s,T%
\right] ;\mathbb{R}\right) \times \mathbb{L}_{\mathcal{F}}^{2}\left( \left[
s,T\right] ;\mathbb{R}\right) \times \mathbb{M}_{\mathcal{F}}^{2}\left( %
\left[ s,T\right] ;\mathbb{R}\right) .$ Moreover when the jump coefficient $%
g\equiv 0$\ the above equations (\ref{2.5})-(\ref{2.6}) reduces to (Buckdahn
et al., \cite{buckdahn1}\ equations (2.7) and (2.10)).$\smallskip $

\noindent Since the derivatives $f_{x},$ $f_{xx},$ $f_{y},$ $\sigma _{x},$ $%
\sigma _{xx},$ $\sigma _{y},$ $\ell _{x},$ $\ell _{y},$ $g_{x},$ $g_{xx}$, $%
h_{x}$, and $h_{y}$ are bounded, by assumptions (H1) and (H2), we have the
following estimate%
\begin{equation}
\mathbb{E}\left[ \sup_{s\leq t\leq T}\left\vert \Psi (t)\right\vert
^{2}+\int_{s}^{T}\left\vert K(t)\right\vert ^{2}dt+\int_{s}^{T}\int_{\Theta
}\left\vert \mathbf{\gamma }_{t}(\theta )\right\vert ^{2}\mu (d\theta )dt%
\right] \leq C.  \tag{2.7}  \label{2.7}
\end{equation}%
\begin{equation}
\mathbb{E}\left[ \sup_{s\leq t\leq T}\left\vert Q(t)\right\vert
^{2}+\int_{s}^{T}\left\vert R(t)\right\vert ^{2}dt+\int_{s}^{T}\int_{\Theta
}\left\vert \Gamma _{t}(\theta )\right\vert ^{2}\mu (d\theta )dt\right] \leq
C.  \tag{2.8}  \label{2.8}
\end{equation}%
\noindent Related with $\left( \Psi ^{\ast }(t),K^{\ast }(t),\gamma
_{t}^{\ast }(\theta )\right) $ we denote%
\begin{equation}
\begin{array}{c}
\delta H(t)=\Psi ^{\ast }(t)\delta f(t)+K^{\ast }(t)\delta \sigma
(t)+\dint_{\Theta }\delta g\left( t,\theta \right) \gamma _{t}^{\ast
}(\theta )\mu (d\theta )-\delta \ell (t),\smallskip \smallskip \\ 
H_{x}(t)=f_{x}\left( t\right) \Psi ^{\ast }(t)+\sigma _{x}\left( t\right)
K^{\ast }(t)+\dint_{\Theta }g_{x}\left( t,\theta \right) \gamma _{t}^{\ast
}(\theta )\mu (d\theta )-\ell _{x}\left( t\right) ,\smallskip \smallskip \\ 
H_{xx}(t)=f_{xx}\left( t\right) \Psi ^{\ast }(t)+\sigma _{xx}\left( t\right)
K^{\ast }(t)+\dint_{\Theta }g_{xx}\left( t,\theta \right) \gamma _{t}^{\ast
}(\theta )\mu (d\theta )-\ell _{xx}\left( t\right) ,%
\end{array}
\tag{2.9}  \label{2.9}
\end{equation}

\section{Stochastic Maximum Principle for Optimality}

\noindent In this section, we obtain a necessary conditions of optimality,
where the system is\ described by nonlinear controlled SDEs of Mean-field
type with jump processes, using spike variation method. The control domain
need not be convex. The proof follows the general ideas as in Buckdahn et
al., \cite{buckdahn1} and Tang et al.,\ \cite{tang}. Note that in \cite%
{buckdahn1} the authors studied the Brownian case only.

\noindent The main result of this paper is stated in the following theorem.

\noindent Let $x^{\ast }(\cdot )$ be the trajectory of the control system (%
\ref{1.1}) corresponding to the optimal control $u^{\ast }(\cdot ),$ and $%
\left( \Psi ^{\ast }(\cdot ),K^{\ast }(\cdot ),\mathbf{\gamma }^{\ast
}(\cdot )\right) ,$\ $\left( Q^{\ast }(\cdot ),R^{\ast }(\cdot ),\Gamma
^{\ast }(\cdot )\right) $\ be the solution of adjoint equations (\ref{2.5})
and (\ref{2.6}) respectively, corresponding to $u^{\ast }(\cdot )$.$%
\smallskip $

\noindent \textbf{Theorem 3.1.}\textit{\ (Stochastic Maximum Principle for
Optimality). }Let Hypotheses (H1) and (H2) hold. If $\left( u^{\ast }(\cdot
),x^{\ast }(\cdot )\right) $\ is an optimal solution of the control problem (%
\ref{1.1})-(\ref{1.2}). Then there are two trible of $\mathcal{F}_{t}%
\mathcal{-}$adapted processes $\left( \Psi ^{\ast }(\cdot ),K^{\ast }(\cdot
),\mathbf{\gamma }^{\ast }(\cdot )\right) $ and $\left( Q^{\ast }(\cdot
),R^{\ast }(\cdot ),\Gamma ^{\ast }(\cdot )\right) $ that satisfy \textit{(%
\ref{2.5}) and (\ref{2.6}) }respectively, such that for all $u\in \mathbb{A}%
: $

\begin{equation}
\begin{array}{l}
H(t,x^{\ast }(t),\mathbb{E}(x^{\ast }(t)\mathbb{)},u,\Psi ^{\ast
}(t),K^{\ast }(t),\gamma _{t}^{\ast }(\theta ))\smallskip \smallskip \\ 
-H(t,x^{\ast }(t),\mathbb{E}(x^{\ast }(t)\mathbb{)},u^{\ast }(t),\Psi ^{\ast
}(t),K^{\ast }(t),\gamma _{t}^{\ast }(\theta ))\smallskip \smallskip \\ 
+\frac{1}{2}\left( \sigma \left( t,x^{\ast }(t),\mathbb{E}(x^{\ast }(t)%
\mathbb{)},u\right) -\sigma \left( t,x^{\ast }(t),\mathbb{E}(x^{\ast }(t)%
\mathbb{)},u^{\ast }(t)\right) \right) ^{2}Q^{\ast }(t)\smallskip \smallskip
\\ 
+\frac{1}{2}\dint_{\Theta }(g\left( t,x^{\ast }(t),u,\theta \right) -g\left(
t,x^{\ast }(t),u^{\ast }(t),\theta \right) )^{2}\left( Q^{\ast }(t)+\Gamma
_{t}^{\ast }(\theta )\right) \mu (d\theta )\leq 0.\smallskip \smallskip \\ 
\text{ \ \ \ \ \ \ \ \ \ \ \ \ \ \ \ \ \ \ \ \ \ \ \ \ \ \ \ \ \ \ \ \ \ \ \
\ \ \ \ \ \ \ \ \ \ \ \ \ \ \ \ \ \ \ \ \ \ \ \ \ \ \ \ }\mathbb{P-}a.s.,%
\text{\ }a.e.\text{ }t\in \left[ s,T\right] .%
\end{array}
\tag{3.1}  \label{3.1}
\end{equation}

\noindent To prove \textit{Theorem 3.1 }we need some preliminary results
given in the following Lemmas.

Let $\left( u^{\ast }(\cdot ),x^{\ast }(\cdot )\right) $\textit{\ }be the
optimal solution of the control problem (\ref{1.1})-(\ref{1.2})\textit{. }%
Following Tang et al., \cite{tang}, and Buckdahn \cite{buckdahn1}, we derive
the variational inequality \textit{(\ref{3.1})} in several steps, from the
fact that%
\begin{equation}
J^{^{s,\zeta }}\left( u^{\varepsilon }(\cdot )\right) -J^{^{s,\zeta }}\left(
u^{\ast }(\cdot )\right) \geq 0,  \tag{3.2}  \label{3.2}
\end{equation}%
where $u^{\varepsilon }(\cdot )$ is the so called spike variation of $%
u^{\ast }(\cdot )$ defined as follows.

\noindent For $\varepsilon >0$, we choose a Borel measurable set $\mathcal{E}%
_{\varepsilon }\subset \left[ s,T\right] $ such that ${\small \upsilon }(%
\mathcal{E}_{\varepsilon })=\varepsilon $, where ${\small \upsilon }(%
\mathcal{E}_{\varepsilon })$ denote the Lebesgue measure of the subset $%
\mathcal{E}_{\varepsilon },$ and we consider the control process which is
the spike variation of $u^{\ast }(\cdot )$%
\begin{equation}
u^{\varepsilon }(t)=\left\{ 
\begin{array}{l}
u:t\in \mathcal{E}_{\varepsilon },\smallskip \smallskip \\ 
u^{\ast }(t):t\in \left[ s,T\right] \mid \mathcal{E}_{\varepsilon },%
\end{array}%
\right.  \tag{3.3}  \label{3.3}
\end{equation}%
where $\varepsilon >0$ is sufficiently small and $u$ is an arbitrary element 
$\mathcal{F}_{t}-$measurable random variable with values in $\mathbb{A}$,
such that $\sup_{w\in \Omega }\left\vert u(w)\right\vert <\infty $, which we
consider as fixed from now on.

\noindent Let $x_{1}^{\varepsilon }(\cdot )$ and $x_{2}^{\varepsilon }(\cdot
)$ be the solutions of the following SDEs respectively%
\begin{equation}
\left\{ 
\begin{array}{l}
dx_{1}^{\varepsilon }(t)=\left\{ f_{x}(t)x_{1}^{\varepsilon }(t)+f_{y}(t)%
\mathbb{E}\left( x_{1}^{\varepsilon }(t)\right) +\delta f(t)\mathbf{I}_{%
\mathcal{E}_{\varepsilon }}(t)\right\} dt\smallskip \smallskip \\ 
\text{ \ \ \ \ \ \ \ \ }+\left\{ \sigma _{x}(t)x_{1}^{\varepsilon
}(t)+\sigma _{y}(t)\mathbb{E}\left( x_{1}^{\varepsilon }(t)\right) +\delta
\sigma (t)\mathbf{I}_{\mathcal{E}_{\varepsilon }}(t)\right\} dW(t)\smallskip
\smallskip \\ 
\text{ \ \ \ \ \ \ \ \ }+\dint_{\Theta }\left\{ g_{x}\left( t_{-},\theta
\right) x_{1}^{\varepsilon }(t)+\delta g(t_{-},\theta )\mathbf{I}_{\mathcal{E%
}_{\varepsilon }}(t)\right\} N\left( d\theta ,dt\right) ,\smallskip
\smallskip \\ 
x_{1}^{\varepsilon }(s)=0,%
\end{array}%
\right.  \tag{3.4}  \label{3.4}
\end{equation}%
and%
\begin{equation}
\left\{ 
\begin{array}{l}
dx_{2}^{\varepsilon }(t)=\left\{ f_{x}(t)x_{2}^{\varepsilon }(t)+f_{y}(t)%
\mathbb{E}\left( x_{2}^{\varepsilon }(t)\right) +\mathcal{L}%
_{t}(f,x_{1}^{\varepsilon })+\delta f_{x}(t)\mathbf{I}_{\mathcal{E}%
_{\varepsilon }}(t)\right\} dt\smallskip \smallskip \\ 
\text{ \ \ \ \ }+\left\{ \sigma _{x}(t)x_{2}^{\varepsilon }(t)+\sigma _{y}(t)%
\mathbb{E}\left( x_{2}^{\varepsilon }(t)\right) +\mathcal{L}_{t}(\sigma
,x_{1}^{\varepsilon })+\delta \sigma _{x}(t)\mathbf{I}_{\mathcal{E}%
_{\varepsilon }}(t)\right\} dW(t)\smallskip \smallskip \\ 
\text{ \ \ \ \ }+\dint_{\Theta }\left\{ g_{x}\left( t_{-},\theta \right)
x_{2}^{\varepsilon }(t)+\mathcal{L}_{t,\theta }(g,x_{1}^{\varepsilon
})+\delta g_{x}(t_{-},\theta )\mathbf{I}_{\mathcal{E}_{\varepsilon
}}(t)\right\} N\left( d\theta ,dt\right) ,\smallskip \smallskip \\ 
x_{2}^{\varepsilon }(s)=0.%
\end{array}%
\right.  \tag{3.5}  \label{3.5}
\end{equation}

\noindent Noting that equation (\ref{3.4}) is called the first-order
variational equation and equation (\ref{3.5}) is called the second-order
variational equation.$\smallskip \smallskip $

\noindent Our first Lemma below deals with the duality relations between $%
\Psi (t),$ $x_{1}^{\varepsilon }(t)$ and $x_{2}^{\varepsilon }(t)$.

\noindent \textbf{Lemma 3.1. }We have

\begin{equation}
\begin{array}{l}
\mathbb{E}\left( \Psi (T)x_{1}^{\varepsilon }(T)\right) =\mathbb{E}%
\dint_{s}^{T}x_{1}^{\varepsilon }(t)\left[ \left( \ell _{x}(t)+\mathbb{E(}%
\ell _{y}(t)\right) \right] dt\smallskip \smallskip \\ 
\text{ \ \ \ \ \ \ \ \ \ \ \ \ \ \ \ \ \ \ \ \ \ \ }+\mathbb{E}%
\dint_{s}^{T}\left\{ \Psi (t)\delta f(t)+K(t)\delta \sigma (t)\right\} 
\mathbf{I}_{\mathcal{E}_{\varepsilon }}(t)dt\smallskip \smallskip \\ 
\text{ \ \ \ \ \ \ \ \ \ \ \ \ \ \ \ \ \ \ \ \ \ \ }+\mathbb{E}%
\dint_{s}^{T}\dint_{\Theta }\gamma _{t}(\theta )\delta g(t,\theta )\mathbf{I}%
_{\mathcal{E}_{\varepsilon }}(t)\mu \left( d\theta \right) dt,%
\end{array}
\tag{3.6}  \label{3.6}
\end{equation}%
and

\begin{equation}
\begin{array}{l}
\mathbb{E}\left( \Psi (T)x_{2}^{\varepsilon }(T)\right) =\mathbb{E}%
\dint_{s}^{T}x_{2}^{\varepsilon }(t)\left[ \left( \ell _{x}(t)+\mathbb{E(}%
\ell _{y}(t)\right) \right] dt\smallskip \smallskip \\ 
\text{ \ \ \ \ \ \ \ \ \ \ \ \ \ \ \ \ \ \ \ \ \ \ }+\mathbb{E}%
\dint_{s}^{T}\left\{ \Psi (t)\delta f_{x}(t)+K(t)\delta \sigma
_{x}(t)\right\} x_{1}^{\varepsilon }(t)\mathbf{I}_{\mathcal{E}_{\varepsilon
}}(t)dt\smallskip \smallskip \\ 
\text{ \ \ \ \ \ \ \ \ \ \ \ \ \ \ \ \ \ \ \ \ \ \ }+\mathbb{E}%
\dint_{s}^{T}\dint_{\Theta }\gamma _{t}(\theta )\delta g_{x}(t,\theta
)x_{1}^{\varepsilon }(t)\mathbf{I}_{\mathcal{E}_{\varepsilon }}(t)\mu \left(
d\theta \right) dt\smallskip \smallskip \\ 
\text{ \ \ \ \ \ \ \ \ \ \ \ \ \ \ \ \ \ \ \ \ \ \ }+\mathbb{E}%
\dint_{s}^{T}\Psi (t)\mathcal{L}_{t}(f,x_{1}^{\varepsilon })+K(t)\mathcal{L}%
_{t}(\sigma ,x_{1}^{\varepsilon })dt\smallskip \smallskip \\ 
\text{ \ \ \ \ \ \ \ \ \ \ \ \ \ \ \ \ \ \ \ \ \ \ }+\mathbb{E}%
\dint_{s}^{T}\dint_{\Theta }\gamma _{t}(\theta )\mathcal{L}_{t,\theta
}(g,x_{1}^{\varepsilon })\mu \left( d\theta \right) dt.%
\end{array}
\tag{3.7}  \label{3.7}
\end{equation}

\noindent \textbf{Proof. }By applying It\^{o}'s formula for jump processes
(see Lemma A1), then we get

\begin{equation}
\begin{array}{l}
\mathbb{E}\left( \Psi (T)x_{1}^{\varepsilon }(T)\right) =\mathbb{E}%
\dint_{s}^{T}\Psi (t)dx_{1}^{\varepsilon }(t)+\mathbb{E}\dint_{s}^{T}x_{1}^{%
\varepsilon }(t)d\Psi (t)\smallskip \smallskip \\ 
\text{ \ \ \ \ \ \ \ \ \ \ \ \ \ \ \ }+\mathbb{E}\dint_{s}^{T}K(t)\left[
\sigma _{x}(t)x_{1}^{\varepsilon }(t)+\sigma _{y}(t)\mathbb{E}\left(
x_{1}^{\varepsilon }(t)\right) +\delta \sigma (t)\mathbf{I}_{\mathcal{E}%
_{\varepsilon }}(t)\right] dt\smallskip \smallskip \\ 
\text{ \ \ \ \ \ \ \ \ \ \ \ \ \ \ \ }+\mathbb{E}\dint_{s}^{T}\dint_{\Theta
}\gamma _{t}(\theta )\left[ g_{x}\left( t,\theta \right) x_{1}^{\varepsilon
}(t)+\delta g(t,\theta )\mathbf{I}_{\mathcal{E}_{\varepsilon }}(t)\right]
\mu \left( d\theta \right) dt\smallskip \smallskip \\ 
\text{ \ \ \ \ \ \ \ \ \ \ \ \ \ \ \ }=I_{1}^{\varepsilon
}+I_{2}^{\varepsilon }+I_{3}^{\varepsilon }+I_{4}^{\varepsilon }.%
\end{array}
\tag{3.8}  \label{3.8}
\end{equation}%
A simple computation shows that%
\begin{equation}
\begin{array}{l}
I_{1}^{\varepsilon }=\mathbb{E}\dint_{s}^{T}\Psi (t)dx_{1}^{\varepsilon
}(t)\smallskip \smallskip \\ 
\text{ \ \ \ }=\mathbb{E}\dint_{s}^{T}\left\{ \Psi
(t)f_{x}(t)x_{1}^{\varepsilon }(t)+\Psi (t)f_{y}(t)\mathbb{E}\left(
x_{1}^{\varepsilon }(t)\right) +\Psi (t)\delta f(t)\mathbf{I}_{\mathcal{E}%
_{\varepsilon }}(t)\right\} dt,%
\end{array}
\tag{3.9}  \label{3.9}
\end{equation}%
and%
\begin{equation}
\begin{array}{l}
I_{2}^{\varepsilon }=\mathbb{E}\dint_{s}^{T}x_{1}^{\varepsilon }(t)d\Psi
(t)\smallskip \smallskip \\ 
\text{ \ \ \ }=-\mathbb{E}\dint_{s}^{T}\left\{ x_{1}^{\varepsilon
}(t)f_{x}\left( t\right) \Psi (t)+x_{1}^{\varepsilon }(t)\mathbb{E}\left(
f_{y}^{\top }(t)\Psi (t)\right) \right. +x_{1}^{\varepsilon }(t)\sigma
_{x}\left( t\right) K(t)\smallskip \smallskip \\ 
\text{ \ \ \ }+\left. x_{1}^{\varepsilon }(t)\mathbb{E}\left( \sigma
_{y}(t)K(t)\right) +x_{1}^{\varepsilon }(t)\ell _{x}\left( t\right)
+x_{1}^{\varepsilon }(t)\mathbb{E}\left( \ell _{y}(t\right) )\right\}
dt\smallskip \smallskip \\ 
\text{ \ \ \ }-\mathbb{E}\dint_{s}^{T}\dint_{\Theta }x_{1}^{\varepsilon
}(t)g_{x}\left( t,\theta \right) \mathbf{\gamma }_{t}(\theta )\mu (d\theta
)dt.\smallskip \smallskip%
\end{array}
\tag{3.10}  \label{3.10}
\end{equation}%
By standard arguments we get

\begin{equation}
\begin{array}{l}
I_{3}^{\varepsilon }=\mathbb{E}\dint_{s}^{T}K(t)\left[ \sigma
_{x}(t)x_{1}^{\varepsilon }(t)+\sigma _{y}(t)\mathbb{E}\left(
x_{1}^{\varepsilon }(t)\right) +\delta \sigma (t)\mathbf{I}_{\mathcal{E}%
_{\varepsilon }}(t)\right] dt\smallskip \smallskip \\ 
\text{ \ \ }=\mathbb{E}\dint_{s}^{T}K(t)\sigma _{x}(t)x_{1}^{\varepsilon
}(t)dt+\mathbb{E}\dint_{s}^{T}K(t)\sigma _{y}(t)\mathbb{E}\left(
x_{1}^{\varepsilon }(t)\right) dt\smallskip \smallskip \\ 
\text{ \ \ }+\mathbb{E}\dint_{s}^{T}K(t)\delta \sigma (t)\mathbf{I}_{%
\mathcal{E}_{\varepsilon }}(t)dt,%
\end{array}
\tag{3.11}  \label{3.11}
\end{equation}

and%
\begin{equation}
\begin{array}{l}
I_{4}^{\varepsilon }=\mathbb{E}\dint_{s}^{T}\dint_{\Theta }\gamma
_{t}(\theta )\left[ g_{x}\left( t,\theta \right) x_{1}^{\varepsilon
}(t)+\delta g(t,\theta )\mathbf{I}_{\mathcal{E}_{\varepsilon }}(t)\right]
\mu \left( d\theta \right) dt\smallskip \smallskip \\ 
\text{ \ \ }=\mathbb{E}\dint_{s}^{T}\dint_{\Theta }\gamma _{t}(\theta
)g_{x}\left( t,\theta \right) x_{1}^{\varepsilon }(t)\mu \left( d\theta
\right) dt\smallskip \smallskip \\ 
\text{ \ \ }+\mathbb{E}\dint_{s}^{T}\dint_{\Theta }\gamma _{t}(\theta
)\delta g(t,\theta )\mathbf{I}_{\mathcal{E}_{\varepsilon }}(t)\mu \left(
d\theta \right) dt.%
\end{array}
\tag{3.12}  \label{3.12}
\end{equation}%
Finally the duality relation (\ref{3.6})\textit{\ }follows from combining (%
\ref{3.9})$\sim $(\ref{3.12}) and (\ref{3.8}). Similarly we can prove second
duality relation (\ref{3.7}).

\noindent To this end we need the following estimations. Let $x^{\varepsilon
}(\cdot )$ be the solutions of the SDEs-(\ref{1.1}) corresponding to the
control $u^{\varepsilon }(\cdot ).$

\noindent \textbf{Lemma 3.2. }Let Hypotheses (H1) and (H2) hold. Then we
have for any $k\geq 1:$

\begin{equation}
\mathbb{E(}\sup_{s\leq t\leq T}\left\vert x_{1}^{\varepsilon }(t)\right\vert
^{2k})\leq C\varepsilon ^{k}.  \tag{3.13}  \label{3.13}
\end{equation}

\begin{equation}
\sup_{s\leq t\leq T}\left\vert \mathbb{E}\left( x_{1}^{\varepsilon
}(t)\right) \right\vert ^{2}\leq \varepsilon \rho \left( \varepsilon \right)
.  \tag{3.14}  \label{3.14}
\end{equation}%
\begin{equation}
\mathbb{E(}\sup_{s\leq t\leq T}\left\vert x_{2}^{\varepsilon }(t)\right\vert
^{2k})\leq C\varepsilon ^{2k}.  \tag{3.15}  \label{3.15}
\end{equation}

\begin{equation}
\mathbb{E(}\sup_{s\leq t\leq T}\left\vert x^{\varepsilon }(t)-x^{\ast
}(t)\right\vert ^{2k})\leq C\varepsilon ^{k}.  \tag{3.16}  \label{3.16}
\end{equation}%
\begin{equation}
\mathbb{E(}\sup_{s\leq t\leq T}\left\vert x^{\varepsilon }(t)-x^{\ast
}(t)-x_{1}^{\varepsilon }(t)\right\vert ^{2k})\leq C\varepsilon ^{2k}. 
\tag{3.17}  \label{3.17}
\end{equation}%
\begin{equation}
\mathbb{E(}\sup_{s\leq t\leq T}\left\vert x^{\varepsilon }(t)-x^{\ast
}(t)-x_{1}^{\varepsilon }(t)-x_{2}^{\varepsilon }(t)\right\vert ^{2k})\leq
C_{k,\mu (\Theta )}\varepsilon ^{2k}\rho _{k}\left( \varepsilon \right) , 
\tag{3.18}  \label{3.18}
\end{equation}%
where $C_{k}$ is a positive constant depend to $k$ and $\rho ,\rho
_{k}:\left( 0,\infty \right) \rightarrow \left( 0,\infty \right) $ such that 
$\rho \left( \varepsilon \right) \rightarrow 0$ and $\rho _{k}\left(
\varepsilon \right) \rightarrow 0$ as $\varepsilon \rightarrow 0.$

\noindent To prove Lemma 3.2 we need some results given in the following
Lemma.$\smallskip $

\noindent \textbf{Lemma 3.3. }For any progressively measurable process $%
\left( \Phi \left( t\right) \right) _{t\in \left[ s,T\right] }$ for which
for any $p>1$, there exists a positive constant $C_{p}$ such that%
\begin{equation}
\mathbb{E(}\sup_{s\leq t\leq T}\left\vert \Phi (t)\right\vert ^{p})\leq
C_{p}.  \tag{3.19}  \label{3.19}
\end{equation}%
Then there exists a function $\widetilde{\rho }:\left( 0,\infty \right)
\rightarrow \left( 0,\infty \right) $ satisfying $\widetilde{\rho }\left(
\varepsilon \right) \rightarrow 0$ as $\varepsilon \rightarrow 0$ such that
for $\varepsilon >0:$

\begin{equation}
\left\vert \mathbb{E}\left( \Phi (T)x_{1}^{\varepsilon }(T)\right)
\right\vert ^{2}+\int_{s}^{T}\left\vert \mathbb{E}\left( \Phi
(t)x_{1}^{\varepsilon }(t)\right) \right\vert ^{2}dt\leq C_{(T,\mu (\Theta
))}\varepsilon \widetilde{\rho }\left( \varepsilon \right) .  \tag{3.20}
\label{3.20}
\end{equation}%
\noindent \textbf{Proof. }First we set for $t\in \left[ s,T\right] :\eta
\left( t\right) =\exp \left\{ Z(t)\right\} ,$ where%
\begin{eqnarray*}
Z\left( t\right) &=&-\int_{s}^{t}\left[ f_{x}(r)-\frac{1}{2}\left\vert
\sigma _{x}(r)\right\vert ^{2}-\frac{1}{2}\int_{\Theta }\left(
g_{x}(r,\theta )\right) ^{2}\mu (d\theta )\right] dr-\int_{s}^{t}\sigma
_{x}(r)dw(r) \\
&&-\int_{s}^{t}\int_{\Theta }g_{x}(r_{-},\theta )N(d\theta ,dr),
\end{eqnarray*}%
and we denote by $\rho \left( t\right) =\eta \left( t\right) ^{-1}=\exp
\left\{ -Z(t)\right\} .$

\noindent By using It\^{o} formula for the exponential $\exp \left\{
Z(t)\right\} $ we get%
\begin{equation*}
d\left( \exp \left\{ Z(t)\right\} \right) =\exp \left\{ Z(t)\right\} dZ(t)+%
\frac{1}{2}\exp \left\{ Z(t)\right\} d\left\langle Z(t);Z(t)\right\rangle ,
\end{equation*}%
this shows that

\begin{equation}
\begin{array}{l}
d\eta \left( t\right) =d\left( \exp \left\{ Z(t)\right\} \right) \smallskip
\smallskip \\ 
\text{ \ \ \ \ \ \ \ \ \ }=-\eta \left( t\right) \left\{ \left[
f_{x}(t)-\left( \sigma _{x}(t)\right) ^{2}-\dint_{\Theta }\left(
g_{x}(t,\theta )\right) ^{2}\mu (d\theta )\right] dt\right. \smallskip
\smallskip \\ 
\text{ \ \ \ \ \ \ \ \ \ \ }+\left. \sigma _{x}(t)dW(t)+\dint_{\Theta
}g_{x}(t_{-},\theta )N(d\theta ,dt)\right\} .%
\end{array}
\tag{3.21}  \label{3.21}
\end{equation}%
\noindent By applying Integration by parts formula for jumps processes $\eta
\left( t\right) x_{1}^{\varepsilon }(t)$ we have%
\begin{eqnarray*}
d\left( \eta \left( t\right) x_{1}^{\varepsilon }(t)\right) &=&\eta \left(
t\right) dx_{1}^{\varepsilon }(t)+x_{1}^{\varepsilon }(t)d\eta \left(
t\right) +d\left\langle \eta \left( t\right) ,x_{1}^{\varepsilon
}(t)\right\rangle , \\
&=&\mathcal{I}_{1}^{\varepsilon }+\mathcal{I}_{2}^{\mathcal{\varepsilon }}+%
\mathcal{I}_{3}^{\varepsilon }.
\end{eqnarray*}%
From (\ref{3.4}) we get%
\begin{eqnarray*}
\mathcal{I}_{1}^{\varepsilon } &=&\eta \left( t\right) dx_{1}^{\varepsilon
}(t) \\
&=&\eta \left( t\right) \left\{ \left[ f_{x}(t)x_{1}^{\varepsilon
}(t)+f_{y}(t)\mathbb{E}\left( x_{1}^{\varepsilon }(t)\right) +\delta f(t)%
\mathbf{I}_{\mathcal{E}_{\varepsilon }}(t)\right] dt\right. \\
&&+\left[ \sigma _{x}(t)x_{1}^{\varepsilon }(t)+\sigma _{y}(t)\mathbb{E}%
\left( x_{1}^{\varepsilon }(t)\right) +\delta \sigma (t)\mathbf{I}_{\mathcal{%
E}_{\varepsilon }}(t)\right] dW(t) \\
&&+\left. \int_{\Theta }\left\{ g_{x}\left( t_{-},\theta \right)
x_{1}^{\varepsilon }(t)+\delta g(t_{-},\theta )\mathbf{I}_{\mathcal{E}%
_{\varepsilon }}(t)\right\} N\left( d\theta ,dt\right) \right\} .
\end{eqnarray*}%
By using (\ref{3.21}) we obtain 
\begin{eqnarray*}
\mathcal{I}_{2}^{\varepsilon } &=&x_{1}^{\varepsilon }(t)d\eta \left(
t\right) \\
&=&-\eta \left( t\right) f_{x}(t)x_{1}^{\varepsilon }(t)dt-\eta \left(
t\right) \sigma _{x}(t)x_{1}^{\varepsilon }(t)dW(t)-\eta \left( t\right)
x_{1}^{\varepsilon }(t)\int_{\Theta }g_{x}(t_{-},\theta )N(d\theta ,dt) \\
&&+\eta \left( t\right) x_{1}^{\varepsilon }(t)\left( \sigma _{x}(t)\right)
^{2}dt+\eta \left( t\right) \int_{\Theta }\left( g_{x}(t_{-},\theta )\right)
^{2}x_{1}^{\varepsilon }(t)\mu (d\theta ),
\end{eqnarray*}%
and a simple computation we get%
\begin{eqnarray*}
\mathcal{I}_{3}^{\varepsilon } &=&d\left\langle \eta \left( t\right)
,x_{1}^{\varepsilon }(t)\right\rangle =-\eta \left( t\right) \sigma _{x}(t)%
\left[ \sigma _{x}(t)x_{1}^{\varepsilon }(t)+\sigma _{y}(t)\mathbb{E}\left(
x_{1}^{\varepsilon }(t)\right) +\delta \sigma (t)\mathbf{I}_{\mathcal{E}%
_{\varepsilon }}(t)\right] dt \\
&&-\int_{\Theta }\eta \left( t\right) g_{x}(t,\theta )\left\{ g_{x}\left(
t,\theta \right) x_{1}^{\varepsilon }(t)+\delta g(t,\theta )\mathbf{I}_{%
\mathcal{E}_{\varepsilon }}(t)\right\} \mu (d\theta )dt.
\end{eqnarray*}%
Consequently, from the above equations we deduce that

\begin{eqnarray*}
d\left( \eta \left( t\right) x_{1}^{\varepsilon }(t)\right) &=&\mathcal{I}%
_{1}+\mathcal{I}_{2}+\mathcal{I}_{3} \\
&=&\eta \left( t\right) \left\{ \left[ f_{y}(t)\mathbb{E}\left(
x_{1}^{\varepsilon }(t)\right) +\delta f(t)\mathbf{I}_{\mathcal{E}%
_{\varepsilon }}(t)\right] dt\right. \\
&&+\left[ \sigma _{y}(t)\mathbb{E}\left( x_{1}^{\varepsilon }(t)\right)
+\delta \sigma (t)\mathbf{I}_{\mathcal{E}_{\varepsilon }}(t)\right] dW(t) \\
&&\left. +\int_{\Theta }\left\{ \delta g(t_{-},\theta )\mathbf{I}_{\mathcal{E%
}_{\varepsilon }}(t)\right\} N\left( d\theta ,dt\right) \right\} \\
&&-\eta \left( t\right) \left\{ \sigma _{x}(t)\left[ \sigma _{y}(t)\mathbb{E}%
\left( x_{1}^{\varepsilon }(t)\right) +\delta \sigma (t)\mathbf{I}_{\mathcal{%
E}_{\varepsilon }}(t)\right] \right. \\
&&+\left. \int_{\Theta }g_{x}(t,\theta )\delta g(t,\theta )\mathbf{I}_{%
\mathcal{E}_{\varepsilon }}(t)\mu (d\theta )\right\} dt,
\end{eqnarray*}%
by integrating the above equation and the fact $\rho \left( t\right) =\eta
\left( t\right) ^{-1}$ we obtain%
\begin{equation}
\begin{array}{l}
x_{1}^{\varepsilon }(t)=\rho \left( t\right) \dint_{s}^{t}\eta \left(
r\right) \left\{ \left[ f_{y}(r)\mathbb{E}\left( x_{1}^{\varepsilon
}(r)\right) +\delta f(t)\mathbf{I}_{\mathcal{E}_{\varepsilon }}(r)\right]
\right. \smallskip \smallskip \\ 
\text{ \ \ \ \ \ \ \ \ }-\sigma _{x}(r)\sigma _{y}(r)\mathbb{E}\left(
x_{1}^{\varepsilon }(r)\right) +\sigma _{x}(r)\delta \sigma (r)\mathbf{I}_{%
\mathcal{E}_{\varepsilon }}(r)\smallskip \smallskip \\ 
\text{ \ \ \ \ \ \ \ \ }-\left. \dint_{\Theta }g_{x}(r,\theta )\delta
g(r,\theta )\mathbf{I}_{\mathcal{E}_{\varepsilon }}(r)\mu (d\theta )\right\}
dr\smallskip \smallskip \\ 
\text{ \ \ \ \ \ \ \ \ }+\rho \left( t\right) \dint_{s}^{t}\eta \left(
r\right) \left[ \sigma _{y}(r)\mathbb{E}\left( x_{1}^{\varepsilon
}(r)\right) +\delta \sigma (r)\mathbf{I}_{\mathcal{E}_{\varepsilon }}(r)%
\right] dW(r)\smallskip \smallskip \\ 
\text{ \ \ \ \ \ \ \ \ }+\rho \left( t\right) \dint_{s}^{t}\dint_{\Theta
}\eta \left( r\right) \delta g(r_{-},\theta )\mathbf{I}_{\mathcal{E}%
_{\varepsilon }}(r)N\left( d\theta ,dr\right) .%
\end{array}
\tag{3.22}  \label{3.22}
\end{equation}%
Since $f_{x},$ $\sigma _{x},$ $g_{x}(\cdot ,\theta )$ are bounded, then by
using (Proposition A1, Appendix) we get: for all $p>1$ there exists a
positive constant $C=C_{(T,p,\mu (\Theta ))}$ such that%
\begin{equation*}
\mathbb{E}\left[ \sup_{s\leq t\leq T}\left\vert \int_{s}^{t}\int_{\Theta
}g_{x}\left( r_{-},\theta \right) N(d\theta ,dr)\right\vert ^{p}\right] \leq
C_{(T,p,\mu (\Theta ))}\mathbb{E}\left[ \int_{s}^{T}\int_{\Theta }\left\vert
g_{x}\left( r,\theta \right) \right\vert ^{p}\mu (d\theta )dr\right] ,
\end{equation*}%
which shows that 
\begin{equation}
\mathbb{E}\left[ \sup_{s\leq t\leq T}(\left\vert \eta \left( t\right)
\right\vert ^{p}+\left\vert \rho \left( t\right) \right\vert ^{p})\right]
\leq C_{(T,p,\mu (\Theta ))}.  \tag{3.23}  \label{3.23}
\end{equation}%
Moreover, it follows from (\ref{3.19}) that%
\begin{equation}
\mathbb{E}\left[ \sup_{s\leq t\leq T}\left\vert \Phi \left( t\right) \rho
\left( t\right) \right\vert ^{p})\right] \leq C_{(T,p,\mu (\Theta ))}. 
\tag{3.24}  \label{3.24}
\end{equation}%
Next, since $\mathcal{F}_{t}=(\mathcal{F}_{t}^{(W,N)})_{t\in \left[ s,T%
\right] }$ then by applying \textit{Martingale Representation Theorem} for
jump processes (see Lemma A2), there exists a unique $\gamma _{t}\left(
\cdot \right) \in \mathbb{L}_{\mathcal{F}}^{2}\left( \left[ s,t\right]
\right) $ and unique $\xi _{t}(\cdot ,\theta )\in \mathbb{M}_{\mathcal{F}%
}^{2}\left( \left[ s,t\right] \right) $ such that $\forall t\in \left[ s,T%
\right] :$%
\begin{equation}
\Phi \left( t\right) \rho \left( t\right) =\mathbb{E}\left( \Phi \left(
t\right) \rho \left( t\right) \right) +\int_{s}^{t}\gamma _{t}\left(
r\right) dW(r)+\int_{s}^{t}\int_{\Theta }\xi _{t}(r,\theta )N(d\theta ,dr).%
\text{ }\mathbb{P}-a.s.  \tag{3.25}  \label{3.25}
\end{equation}%
Noting that, for every $p>1,$ with the help of (\ref{3.22}) it follows from
the Bulkholder-Davis-Gundy inequality and Proposition A1 that there exists a
constant $C_{(T,p,\mu (\Theta ))}$ such that: for $p>1,$%
\begin{eqnarray*}
&&\mathbb{E}\left[ \left( \int_{s}^{t}\left\vert \gamma _{t}\left( r\right)
\right\vert ^{2}dr\right) ^{\frac{p}{2}}\right] +\mathbb{E}\left[ \left(
\int_{s}^{t}\int_{\Theta }\left\vert \xi _{t}(r,\theta )\right\vert ^{2}\mu
(d\theta )dr\right) ^{\frac{p}{2}}\right] \\
&\leq &C_{p}\mathbb{E}\left[ \sup_{s\leq \tau \leq t}\left\vert
\int_{s}^{\tau }\gamma _{t}\left( r\right) dW(r)\right\vert ^{p}\right] \\
&&+C_{(T,p,\mu (\Theta ))}\mathbb{E}\left[ \sup_{s\leq \tau \leq
t}\left\vert \int_{s}^{\tau }\int_{\Theta }\xi _{t}(r_{-},\theta )N(d\theta
,dr)\right\vert ^{p}\right] \\
&\leq &C_{p}\left( 1+\tfrac{1}{p-1}\right) ^{p}\mathbb{E}\left[ \left\vert
\int_{s}^{t}\gamma _{t}\left( r\right) dW(r)\right\vert ^{p}\right] \\
&&+C_{(T,p,\mu (\Theta ))}\mathbb{E}\left[ \left\vert
\int_{s}^{t}\int_{\Theta }\xi _{t}(r_{-},\theta )N(d\theta ,dr)\right\vert
^{p}\right] \\
&\leq &C_{(T,p,\mu (\Theta ))}\mathbb{E}\left[ \left\vert \Phi \left(
t\right) \rho \left( t\right) -\mathbb{E}\left( \Phi \left( t\right) \rho
\left( t\right) \right) \right\vert ^{p}\right] \\
&\leq &C_{(T,p,\mu (\Theta ))}\left\{ \mathbb{E(}\left\vert \Phi \left(
t\right) \rho \left( t\right) \right\vert ^{p})+\left\vert \mathbb{E}\left(
\Phi \left( t\right) \rho \left( t\right) \right) \right\vert ^{p}\right\} \\
&\leq &C_{(T,p,\mu (\Theta ))}\mathbb{E}\left[ \left\vert \Phi \left(
t\right) \rho \left( t\right) \right\vert ^{p}\right] \\
&\leq &C_{(T,p,\mu (\Theta ))}\mathbb{E}\left[ \sup_{s\leq t\leq
T}\left\vert \Phi \left( t\right) \rho \left( t\right) \right\vert ^{p}%
\right] \leq C_{(T,p,\mu (\Theta ))}.
\end{eqnarray*}%
This shows that%
\begin{equation}
\sup_{s\leq t\leq T}\mathbb{E}\left[ \left( \int_{s}^{t}\left\vert \gamma
_{t}\left( r\right) \right\vert ^{2}dr\right) ^{\frac{p}{2}}\right] \leq
C_{(T,p,\mu (\Theta ))},  \tag{3.26}  \label{3.26}
\end{equation}%
and%
\begin{equation}
\sup_{s\leq t\leq T}\mathbb{E}\left[ \left( \int_{s}^{t}\int_{\Theta
}\left\vert \xi _{t}(r,\theta )\right\vert ^{2}\mu (d\theta )dr\right) ^{%
\frac{p}{2}}\right] \leq C_{(T,p,\mu (\Theta ))}.  \tag{3.27}  \label{3.27}
\end{equation}%
\noindent Now we consider%
\begin{equation}
\Phi \left( t\right) x_{1}^{\varepsilon }\left( t\right) =\mathcal{J}%
_{1}^{\varepsilon }(t)+\mathcal{J}_{2}^{\varepsilon }(t)+\mathcal{J}%
_{3}^{\varepsilon }(t),\text{ }t\in \left[ s,T\right] ,  \tag{3.28}
\label{3.28}
\end{equation}%
where%
\begin{eqnarray*}
\mathcal{J}_{1}^{\varepsilon }(t) &=&\Phi \left( t\right) \rho \left(
t\right) \int_{s}^{t}\eta \left( r\right) \left\{ \left[ f_{y}(r)\mathbb{E}%
\left( x_{1}^{\varepsilon }(r)\right) +\delta f(t)\mathbf{I}_{\mathcal{E}%
_{\varepsilon }}(r)\right] \right. \\
&&-\sigma _{x}(r)\sigma _{y}(r)\mathbb{E}\left( x_{1}^{\varepsilon
}(r)\right) +\sigma _{x}(r)\delta \sigma (r)\mathbf{I}_{\mathcal{E}%
_{\varepsilon }}(r) \\
&&-\left. \int_{\Theta }g_{x}(r_{-},\theta )\delta g(r,\theta )\mathbf{I}_{%
\mathcal{E}_{\varepsilon }}(r)\mu (d\theta )\right\} dr,
\end{eqnarray*}%
\begin{equation*}
\mathcal{J}_{2}^{\varepsilon }(t)=\Phi \left( t\right) \rho \left( t\right)
\int_{s}^{t}\eta \left( r\right) \left[ \sigma _{y}(r)\mathbb{E}\left(
x_{1}^{\varepsilon }(r)\right) +\delta \sigma (r)\mathbf{I}_{\mathcal{E}%
_{\varepsilon }}(r)\right] dW(r),
\end{equation*}%
and%
\begin{equation*}
\mathcal{J}_{3}^{\varepsilon }(t)=\Phi \left( t\right) \rho \left( t\right)
\int_{s}^{t}\int_{\Theta }\eta \left( r\right) \delta g(r,\theta )\mathbf{I}%
_{\mathcal{E}_{\varepsilon }}(r)N\left( d\theta ,dr\right) .
\end{equation*}%
We estimate now the first term in the right-hand side of (\ref{3.28}).
First, since $f_{y},$ $\sigma _{x}$, $\sigma _{y}$, are bounded and fact
that $\sup_{\theta \in \Theta }\left\vert g_{x}(t,\theta )\right\vert
<+\infty $ (see H2) we get%
\begin{eqnarray*}
\left\vert \mathbb{E}\left( \mathcal{J}_{1}^{\varepsilon }(t)\right)
\right\vert &=&\left\vert \mathbb{E}\left\{ \Phi \left( t\right) \rho \left(
t\right) \int_{s}^{t}\eta \left( r\right) \left[ \left( f_{y}(r)-\sigma
_{x}(r)\sigma _{y}(r)\right) \mathbb{E}\left( x_{1}^{\varepsilon }(r)\right)
\right. \right. \right. \\
&&+\left( \delta f(t)+\sigma _{x}(r)\delta \sigma (r)\right) \mathbf{I}_{%
\mathcal{E}_{\varepsilon }}(r) \\
&&-\left. \left. \left. \int_{\Theta }g_{x}(r_{-},\theta )\delta g(r,\theta )%
\mathbf{I}_{\mathcal{E}_{\varepsilon }}(r)\mu (d\theta )\right] dr\right\}
\right\vert \\
&\leq &C_{(\mu (\Theta ))}\mathbb{E}\left\{ \sup_{t\in \left[ s,T\right]
}\left\vert \Phi \left( t\right) \rho \left( t\right) \right\vert \sup_{t\in %
\left[ s,T\right] }\left\vert \eta \left( t\right) \right\vert \left[
\int_{s}^{t}\left\vert \mathbb{E}\left( x_{1}^{\varepsilon }(r)\right)
\right\vert dr+\varepsilon \right] \right\} ,
\end{eqnarray*}%
applying Cauchy-Schwarz inequality,\ then from (\ref{3.23}) and (\ref{3.24})
(with $p=2),$ we get%
\begin{eqnarray*}
\left\vert \mathbb{E}\left( \mathcal{J}_{1}^{\varepsilon }(t)\right)
\right\vert &\leq &C_{(\mu (\Theta ))}\left[ \mathbb{E}\left( \sup_{t\in %
\left[ s,T\right] }\left\vert \Phi \left( t\right) \rho \left( t\right)
\right\vert ^{2}\right) \right] ^{\frac{1}{2}} \\
&&\times \left[ \mathbb{E}\left( \sup_{t\in \left[ s,T\right] }\left\vert
\eta \left( t\right) \right\vert ^{2}\right) \right] ^{\frac{1}{2}}\left[
\int_{s}^{t}\left\vert \mathbb{E}\left( x_{1}^{\varepsilon }(r)\right)
\right\vert dr+\varepsilon \right] \\
&\leq &C_{(T,\mu (\Theta ))}\left[ \int_{s}^{t}\left\vert \mathbb{E}\left(
x_{1}^{\varepsilon }(r)\right) \right\vert dr+\varepsilon \right] ,
\end{eqnarray*}%
by applying Cauchy-Schwarz inequality and the fact that $\left( a+b\right)
^{2}\leq 2a^{2}+2b^{2}$ we can shows that%
\begin{equation}
\begin{array}{l}
\left\vert \mathbb{E}\left( \mathcal{J}_{1}^{\varepsilon }(t)\right)
\right\vert ^{2}\leq C_{(T,\mu (\Theta ))}\left[ 2\left(
\dint_{s}^{t}\left\vert \mathbb{E}\left( x_{1}^{\varepsilon }(r)\right)
\right\vert dr\right) ^{2}+2\varepsilon ^{2}\right] \smallskip \smallskip \\ 
\text{ \ \ \ \ \ \ \ \ \ \ \ \ \ \ \ \ \ }\leq C_{(T,\mu (\Theta ))}\left[
\dint_{s}^{t}\left\vert \mathbb{E}\left( x_{1}^{\varepsilon }(r)\right)
\right\vert ^{2}dr+\varepsilon ^{2}\right] .%
\end{array}
\tag{3.29}  \label{3.29}
\end{equation}%
Next, we proceed to estimate the second term $\mathcal{J}_{2}^{\varepsilon
}(t).$ With the help of (\ref{3.25}) and the \emph{It\^{o} Isometry} we can
get%
\begin{eqnarray*}
\mathbb{E}\left( \mathcal{J}_{2}^{\varepsilon }(t)\right) &=&\mathbb{E}%
\left\{ \Phi \left( t\right) \rho \left( t\right) \int_{s}^{t}\eta \left(
r\right) \left[ \sigma _{y}(r)\mathbb{E}\left( x_{1}^{\varepsilon
}(r)\right) +\delta \sigma (r)\mathbf{I}_{\mathcal{E}_{\varepsilon }}(r)%
\right] dW(r)\right\} \\
&=&\mathbb{E}\left\{ \left[ \mathbb{E}\left( \Phi \left( t\right) \rho
\left( t\right) \right) +\int_{s}^{t}\gamma _{t}\left( r\right)
dW(r)+\int_{s}^{t}\int_{\Theta }\xi _{t}(r_{-},\theta )N(d\theta ,dr)\right]
\right. \\
&&\times \left. \int_{s}^{t}\eta \left( r\right) \left[ \sigma _{y}(r)%
\mathbb{E}\left( x_{1}^{\varepsilon }(r)\right) +\delta \sigma (r)\mathbf{I}%
_{\mathcal{E}_{\varepsilon }}(r)\right] dW(r)\right\} \\
&=&\mathbb{E}\left[ \int_{s}^{t}\gamma _{t}\left( r\right) \eta \left(
r\right) \sigma _{y}(r)\mathbb{E}\left( x_{1}^{\varepsilon }(r)\right) dr%
\right] +\mathbb{E}\left[ \int_{s}^{t}\gamma _{t}\left( r\right) \eta \left(
r\right) \delta \sigma (r)\mathbf{I}_{\mathcal{E}_{\varepsilon }}(r)dr\right]
.
\end{eqnarray*}%
We estimate now the first term in the right hand side of the above equality.
Applying (\ref{3.23})-(\ref{3.26}) then we can get immediately

\begin{equation}
\left\vert \mathbb{E}\int_{s}^{t}\gamma _{t}\left( r\right) \eta \left(
r\right) \sigma _{y}(r)\mathbb{E}\left( x_{1}^{\varepsilon }(r)\right)
dr\right\vert ^{2}\leq C\int_{s}^{t}\left\vert \mathbb{E}\left(
x_{1}^{\varepsilon }(r)\right) \right\vert ^{2}dr,  \tag{3.30}  \label{3.30}
\end{equation}%
however, the second term satisfies%
\begin{equation}
\int_{s}^{T}\left\{ \left\vert \mathbb{E}\left[ \int_{s}^{t}\gamma
_{t}\left( r\right) \eta \left( r\right) \delta \sigma (r)\mathbf{I}_{%
\mathcal{E}_{\varepsilon }}(r)dr\right] \right\vert ^{2}\right\} dt\leq
C\varepsilon \rho _{1}\left( \varepsilon \right) ,  \tag{3.31}  \label{3.31}
\end{equation}%
where%
\begin{equation*}
\rho _{1}\left( \varepsilon \right) =\left\{ \mathbb{E}\left[ \left(
\int_{s}^{T}\int_{s}^{t}\left\vert \gamma _{t}\left( r\right) \right\vert
^{2}\mathbf{I}_{\mathcal{E}_{\varepsilon }}(t)drdt\right) ^{2}\right]
\right\} ^{\frac{1}{2}}.
\end{equation*}

\noindent Noting that since $\lim_{\varepsilon \rightarrow 0}\mathbf{I}_{%
\mathcal{E}_{\varepsilon }}(t)=0$ in measure $dtd\mathbb{P}$ then the by 
\textit{Dominate Convergence Theorem} we get that $\lim_{\varepsilon
\rightarrow 0}\rho _{1}\left( \varepsilon \right) =0.$

\noindent Let us turn to estimate the third term $\mathcal{J}%
_{3}^{\varepsilon }(t)$. Then by Cauchy-Schwarz inequality, we obtain%
\begin{eqnarray*}
\left\vert \mathbb{E}\left( \mathcal{J}_{3}^{\varepsilon }(t)\right)
\right\vert ^{2} &=&\left\vert \mathbb{E}\left\{ \Phi \left( t\right) \rho
\left( t\right) \int_{s}^{t}\int_{\Theta }\eta \left( r\right) \delta
g(r_{-},\theta )\mathbf{I}_{\mathcal{E}_{\varepsilon }}(r)N\left( d\theta
,dr\right) \right\} \right\vert ^{2} \\
&\leq &\left\vert \mathbb{E}\left\{ \sup_{t\in \left[ s,T\right] }\left\vert
\Phi \left( t\right) \rho \left( t\right) \right\vert \sup_{t\in \left[ s,T%
\right] }\left\vert \eta \left( t\right) \right\vert
\int_{s}^{t}\int_{\Theta }\delta g(r_{-},\theta )\mathbf{I}_{\mathcal{E}%
_{\varepsilon }}(r)N\left( d\theta ,dr\right) \right\} \right\vert ^{2} \\
&\leq &C\mathbb{E}\left[ \sup_{t\in \left[ s,T\right] }\left\vert \Phi
\left( t\right) \rho \left( t\right) \right\vert ^{2}\right] \left[ \mathbb{E%
}\sup_{t\in \left[ s,T\right] }\left\vert \eta \left( t\right) \right\vert
^{2}\right] \\
&&\times \mathbb{E}\left[ \left\vert \int_{s}^{t}\int_{\Theta }\delta
g(r_{-},\theta )\mathbf{I}_{\mathcal{E}_{\varepsilon }}(r)N\left( d\theta
,dr\right) \right\vert \right] ^{2},
\end{eqnarray*}%
by applying Propositions A1 then from (\ref{3.23}) and (\ref{3.24}) (with $%
p=2),$ we get%
\begin{equation}
\begin{array}{l}
\left\vert \mathbb{E}\left( \mathcal{J}_{3}^{\varepsilon }(t)\right)
\right\vert ^{2}\leq C_{(T,\mu (\Theta ))}\dint_{s}^{t}\dint_{\Theta
}\left\vert \delta g(r,\theta )\mathbf{I}_{\mathcal{E}_{\varepsilon
}}(r)\right\vert ^{2}\mu (d\theta )dr\smallskip \smallskip \\ 
\text{ \ \ \ \ \ \ \ \ \ \ \ \ \ \ \ \ \ \ }\leq C_{(T,\mu (\Theta
))}\dint_{s}^{t}\sup_{\theta \in \Theta }\left\vert \delta g(r,\theta
)\right\vert ^{2}\dint_{\Theta }\mathbf{I}_{\mathcal{E}_{\varepsilon
}}(r)\mu (d\theta )dr\smallskip \smallskip \\ 
\text{ \ \ \ \ \ \ \ \ \ \ \ \ \ \ \ \ \ \ }\leq C_{(T,\mu (\Theta
))}\varepsilon .%
\end{array}
\tag{3.32}  \label{3.32}
\end{equation}%
Combining (\ref{3.30})$\sim $(\ref{3.32}) and the fact that%
\begin{equation*}
\left\vert \mathbb{E}\left( \Phi \left( t\right) x_{1}^{\varepsilon }\left(
t\right) \right) \right\vert ^{2}\leq 2\left\vert \mathbb{E}\left( \mathcal{J%
}_{1}^{\varepsilon }(t)\right) \right\vert ^{2}+4\left\vert \mathbb{E}\left( 
\mathcal{J}_{2}^{\varepsilon }(t)\right) \right\vert ^{2}+4\left\vert 
\mathbb{E}\left( \mathcal{J}_{3}^{\varepsilon }(t\right) )\right\vert
^{2},t\in \left[ s,T\right] ,
\end{equation*}%
we conclude%
\begin{equation}
\begin{array}{c}
\left\vert \mathbb{E}\left( \Phi \left( t\right) x_{1}^{\varepsilon }\left(
t\right) \right) \right\vert ^{2}\leq C_{(T,\mu (\Theta ))}\left[
\varepsilon ^{2}+\varepsilon +\left\vert \mathbb{E}\left[ \int_{s}^{t}\gamma
_{t}\left( r\right) \eta \left( r\right) \delta \sigma (r)\mathbf{I}_{%
\mathcal{E}_{\varepsilon }}(r)dr\right] \right\vert ^{2}\right. \smallskip
\smallskip \\ 
\left. +\int_{s}^{t}\left\vert \mathbb{E}\left( x_{1}^{\varepsilon }\left(
r\right) \right) \right\vert ^{2}dr\right] ,%
\end{array}
\tag{3.33}  \label{3.33}
\end{equation}%
integrating the above inequality, then with the help of (\ref{3.31}) we get 
\begin{equation}
\int_{s}^{t}\left\vert \mathbb{E}\left( \Phi \left( r\right)
x_{1}^{\varepsilon }\left( r\right) \right) \right\vert ^{2}dr\leq C_{(T,\mu
(\Theta ))}\left[ \varepsilon ^{2}+\varepsilon +\varepsilon \rho _{1}\left(
\varepsilon \right) +\int_{s}^{t}\int_{s}^{e}\left\vert \mathbb{E}\left(
x_{1}^{\varepsilon }\left( r\right) \right) \right\vert ^{2}drde\right] . 
\tag{3.34}  \label{3.34}
\end{equation}%
Now, taking $\Phi \left( t\right) =1$ in (\ref{3.34}) and from Gronwall's
Lemma we have%
\begin{equation}
\int_{s}^{t}\left\vert \mathbb{E}\left( x_{1}^{\varepsilon }\left( r\right)
\right) \right\vert ^{2}dr\leq C_{(T,\mu (\Theta ))}\left( \varepsilon
^{2}+\varepsilon +\varepsilon \rho _{1}\left( \varepsilon \right) \right) . 
\tag{3.35}  \label{3.35}
\end{equation}%
Consequently, from (\ref{3.34}) it holds that%
\begin{equation}
\int_{s}^{t}\left\vert \mathbb{E}\left( \Phi \left( r\right)
x_{1}^{\varepsilon }\left( r\right) \right) \right\vert ^{2}dr\leq C_{(T,\mu
(\Theta ))}\left( \varepsilon ^{2}+\varepsilon +\varepsilon \rho _{1}\left(
\varepsilon \right) \right) .  \tag{3.36}  \label{3.36}
\end{equation}%
Furthermore, from (\ref{3.23}), then by simple computation (with $t=T$) we
can shows that

\begin{equation}
\left\vert \mathbb{E}\left[ \int_{s}^{T}\gamma _{t}\left( r\right) \eta
\left( r\right) \delta \sigma (r)\mathbf{I}_{\mathcal{E}_{\varepsilon }}(r)dr%
\right] \right\vert ^{2}\leq C\varepsilon \rho _{T}\left( \varepsilon
\right) ,  \tag{3.37}  \label{3.37}
\end{equation}%
where%
\begin{equation*}
\rho _{T}\left( \varepsilon \right) =\left\{ \mathbb{E}\left[ \left(
\int_{s}^{T}\left\vert \gamma _{T}\left( r\right) \right\vert ^{2}\mathbf{I}%
_{\mathcal{E}_{\varepsilon }}(t)dr\right) ^{2}\right] \right\} ^{\frac{1}{2}%
}.
\end{equation*}%
Noting that since $\lim_{\varepsilon \rightarrow 0}\mathbf{I}_{\mathcal{E}%
_{\varepsilon }}(t)=0$ in measure $dtd\mathbb{P}$ then with the help \textit{%
Dominate Convergence Theorem }we can shows that $\lim_{\varepsilon
\rightarrow 0}\rho _{T}\left( \varepsilon \right) =0.$

\noindent By combining (\ref{3.33}), (\ref{3.35}) and (\ref{3.37}) we
conclude%
\begin{equation}
\left\vert \mathbb{E}\left( \Phi \left( T\right) x_{1}^{\varepsilon }\left(
T\right) \right) \right\vert ^{2}\leq C_{(T,\mu (\Theta ))}\left(
\varepsilon +\varepsilon ^{2}+\varepsilon \rho _{1}\left( \varepsilon
\right) +\varepsilon \rho _{T}\left( \varepsilon \right) \right) . 
\tag{3.38}  \label{3.38}
\end{equation}%
Finally by setting $\widetilde{\rho }\left( \varepsilon \right) =\left(
\varepsilon +\varepsilon ^{2}+\varepsilon \rho _{1}\left( \varepsilon
\right) +\varepsilon \rho _{T}\left( \varepsilon \right) \right) \rightarrow
0,$ $\varepsilon \rightarrow 0$, then the desired result (\ref{3.20})
follows immediately from (\ref{3.36}) and (\ref{3.38}). This completes the
proof of Lemma 3.3.

\noindent \textbf{Proof of Lemma 3.2.}

\noindent \textit{Proof of estimate (\ref{3.14}): }using (\ref{3.4}) it
holds that%
\begin{equation*}
\mathbb{E(}x_{1}^{\varepsilon }\left( t\right) )=\int_{s}^{t}\left\{ \mathbb{%
E}\left[ f_{x}(r)x_{1}^{\varepsilon }(r)\right] +\mathbb{E}\left(
f_{y}(r)\right) \mathbb{E}\left( x_{1}^{\varepsilon }(r)\right) +\mathbb{E}%
\left( \delta f(r)\mathbf{I}_{\mathcal{E}_{\varepsilon }}(r)\right) \right\}
dr,
\end{equation*}%
then we have%
\begin{equation}
\begin{array}{l}
\left\vert \mathbb{E(}x_{1}^{\varepsilon }\left( t\right) )\right\vert
^{2}\leq 2\left\vert \dint_{s}^{t}\mathbb{E}\left[ f_{x}(r)x_{1}^{%
\varepsilon }(r)\right] dr\right\vert ^{2}\smallskip \smallskip \\ 
\text{ \ \ \ \ \ \ \ \ \ \ \ \ \ \ \ \ \ \ }+2\left\vert \dint_{s}^{t}\left( 
\mathbb{E}\left( f_{y}(r)\right) \mathbb{E}\left( x_{1}^{\varepsilon
}(r)\right) +\mathbb{E}\left( \delta f(r)\mathbf{I}_{\mathcal{E}%
_{\varepsilon }}(r)\right) \right) dr\right\vert ^{2},%
\end{array}
\tag{3.39}  \label{3.39}
\end{equation}%
by setting $\Phi \left( t\right) =f_{x}(t)$ in (\ref{3.36}), then by the
helps of Cauchy-Schwarz inequality and fact that $t\leq T$ we get%
\begin{equation*}
\left\vert \int_{s}^{t}\mathbb{E}\left[ f_{x}(r)x_{1}^{\varepsilon }(r)%
\right] \right\vert ^{2}\leq T\int_{s}^{t}\left\vert \mathbb{E}\left[
f_{x}(r)x_{1}^{\varepsilon }(r)\right] \right\vert ^{2}dr\leq C_{(T,\mu
(\Theta ))}\left( \varepsilon ^{2}+\varepsilon +\varepsilon \rho _{1}\left(
\varepsilon \right) \right) ,
\end{equation*}%
thus, in view of assumption (H1), then from (\ref{3.39}) we obtain%
\begin{equation*}
\left\vert \mathbb{E(}x_{1}^{\varepsilon }\left( t\right) )\right\vert \leq %
\left[ C_{(T,\mu (\Theta ))}\left( \varepsilon ^{2}+\varepsilon +\varepsilon
\rho _{1}\left( \varepsilon \right) \right) \right] ^{\frac{1}{2}%
}+C\int_{s}^{t}\left( \varepsilon +\mathbb{E}\left\vert x_{1}^{\varepsilon
}(r)\right\vert \right) dr.
\end{equation*}%
Finally by applying Gronwall's Lemma, the estimate (\ref{3.14}) follows with 
$\rho \left( \varepsilon \right) =C_{(T,s,\mu (\Theta ))}\left(
1+\varepsilon +\rho _{1}\left( \varepsilon \right) \right) .$

\noindent \textit{Proof of estimate (\ref{3.18}):} First we set%
\begin{equation}
\lambda ^{\varepsilon }(t):=x^{\varepsilon }(t)-x^{\ast
}(t)-x_{1}^{\varepsilon }(t)-x_{2}^{\varepsilon }(t).  \tag{3.40}
\label{3.40}
\end{equation}%
From SDEs (\ref{1.1}), (\ref{3.4}) and (\ref{3.5}) we get 
\begin{equation}
d\lambda ^{\varepsilon }(t)=\Pi _{f}^{\varepsilon }\left( t\right) dt+\Pi
_{\sigma }^{\varepsilon }\left( t\right) dW(t)+\int_{\Theta }\Lambda
_{g}^{\varepsilon }\left( t,\theta \right) N(d\theta ,dt),  \tag{3.41}
\label{3.41}
\end{equation}%
where for $\varphi =f,$ $\sigma ,$ $\ell $%
\begin{equation}
\begin{array}{l}
\Pi _{\varphi }^{\varepsilon }\left( t\right) =\varphi (t,x^{\varepsilon
}(t),\mathbb{E}(x^{\varepsilon }(t)),u^{\varepsilon }(t))-\varphi (t,x^{\ast
}(t),\mathbb{E}(x^{\ast }(t)),u^{\ast }(t))\smallskip \smallskip \\ 
\text{ \ \ \ \ \ \ \ \ \ \ }-\varphi _{x}\left( t\right) \left(
x_{1}^{\varepsilon }(t)+x_{2}^{\varepsilon }(t)\right) -\left\{ \varphi
_{y}\left( t\right) \mathbb{E}\left( x_{1}^{\varepsilon
}(t)+x_{2}^{\varepsilon }(t)\right) \right. \smallskip \smallskip \\ 
\text{ \ \ \ \ \ \ \ \ \ \ }\left. +\mathcal{L}_{t}\left( \varphi
,x_{1}^{\varepsilon }\right) +\left( \delta \varphi \left( t\right) +\delta
\varphi _{x}\left( t\right) x_{1}^{\varepsilon }\left( t\right) \right) 
\mathbf{I}_{\mathcal{E}_{\varepsilon }}(t)\right\} ,%
\end{array}
\tag{3.42}  \label{3.42}
\end{equation}%
and 
\begin{equation}
\begin{array}{l}
\Lambda _{g}^{\varepsilon }\left( t,\theta \right) =g\left( t,x^{\varepsilon
}(t_{-}),u^{\varepsilon }(t),\theta \right) -\left\{ g\left( t,x^{\ast
}(t_{-}),u^{\ast }(t),\theta \right) \right. \smallskip \smallskip \\ 
\text{ \ \ \ }+\left. g_{x}\left( t,\theta \right) \left[ x_{1}^{\varepsilon
}(t)+x_{2}^{\varepsilon }(t)\right] +\mathcal{L}_{t,\theta
}(g,x_{1}^{\varepsilon })+\left[ \delta g(t,\theta )+\delta g_{x}(t,\theta )%
\right] \mathbf{I}_{\mathcal{E}_{\varepsilon }}(t)\right\} .%
\end{array}
\tag{3.43}  \label{3.43}
\end{equation}%
First we estimate the term $\Pi _{\varphi }^{\varepsilon }\left( t\right) .$

\noindent \textit{Estimates of }$\Pi _{\varphi }^{\varepsilon }\left(
t\right) :$ 
\begin{eqnarray*}
&&\varphi (t,x^{\varepsilon }(t),\mathbb{E}(x^{\varepsilon
}(t)),u^{\varepsilon }(t))-\varphi (t,x^{\ast }(t),\mathbb{E}(x^{\ast
}(t)),u^{\ast }(t)) \\
&=&\int_{0}^{1}\left[ \varphi _{x}^{e}(t)\left( x^{\varepsilon }(t)-x^{\ast
}(t)\right) +\varphi _{y}^{e}(t)\left( \mathbb{E}\left( x^{\varepsilon
}(t)\right) -\mathbb{E}\left( x^{\ast }(t)\right) \right) \right] de,
\end{eqnarray*}%
where, for the subscript $\varkappa $ which indicates the first and the
second order derivatives of $\varphi $, respectively, with respect to $%
\varkappa =x,$ $xx,$ $y,xy,$ $yy$, and for real\ $\hslash \in \left[ 0,1%
\right] :$%
\begin{equation*}
\varphi _{\varkappa }^{\hslash }(t)=\varphi _{\varkappa }(t,x^{\ast
}(t)+\hslash \left( x^{\varepsilon }(t)-x^{\ast }(t)\right) ,\mathbb{E}%
\left( x^{\ast }(t)+\hslash \left( x^{\varepsilon }(t)-x^{\ast }(t)\right)
\right) ,u^{\varepsilon }(t)).
\end{equation*}%
Moreover,%
\begin{eqnarray*}
&&\varphi (t,x^{\varepsilon }(t),\mathbb{E}(x^{\varepsilon
}(t)),u^{\varepsilon }(t))-\varphi (t,x^{\ast }(t),\mathbb{E}(x^{\ast
}(t)),u^{\ast }(t)) \\
&&-\left[ \varphi _{x}\left( t\right) \left( x_{1}^{\varepsilon
}(t)+x_{2}^{\varepsilon }(t)\right) \right. \left. +\varphi _{y}\left(
t\right) \mathbb{E}\left( x_{1}^{\varepsilon }(t)+x_{2}^{\varepsilon
}(t)\right) \right] \\
&=&\int_{0}^{1}\left\{ \varphi _{x}^{e}(t)\lambda ^{\varepsilon }(t)+\varphi
_{y}^{e}(t)\mathbb{E}\left( \lambda ^{\varepsilon }(t)\right) +\left(
\varphi _{x}^{e}(t)-\varphi _{x}(t)\right) \left( x_{1}^{\varepsilon
}(t)+x_{2}^{\varepsilon }(t)\right) \right. \\
&&+\left. \left( \varphi _{y}^{e}(t)-\varphi _{y}(t)\right) \mathbb{E}\left(
x_{1}^{\varepsilon }(t)+x_{2}^{\varepsilon }(t)\right) \right\} de.
\end{eqnarray*}%
By similar arguments we get%
\begin{eqnarray*}
\varphi _{x}^{e}(t)-\varphi _{x}(t) &=&e\int_{0}^{1}\left\{ \varphi
_{xx}^{e,\alpha }(t)\left( x^{\varepsilon }(t)-x^{\ast }(t)\right) +\varphi
_{xy}^{e,\alpha }(t)\mathbb{E}\left( x^{\varepsilon }(t)-x^{\ast }(t)\right)
\right\} d\alpha \\
&&+\delta \varphi _{x}(t)\mathbf{I}_{\mathcal{E}_{\varepsilon }}(t) \\
&=&e\int_{0}^{1}\left\{ \varphi _{xx}^{e,\alpha }(t)\lambda ^{\varepsilon
}(t)+\varphi _{xy}^{e,\alpha }(t)\mathbb{E}\left( \lambda ^{\varepsilon
}(t)\right) \right\} d\alpha \\
&&+e\int_{0}^{1}\left\{ \varphi _{xx}^{e,\alpha }(t)\left(
x_{1}^{\varepsilon }(t)+x_{2}^{\varepsilon }(t)\right) \right\} d\alpha \\
&&+e\int_{0}^{1}\left\{ \varphi _{xy}^{e,\alpha }(t)\mathbb{E}\left(
x_{1}^{\varepsilon }(t)+x_{2}^{\varepsilon }(t)\right) \right\} d\alpha
+\delta \varphi _{x}(t)\mathbf{I}_{\mathcal{E}_{\varepsilon }}(t),
\end{eqnarray*}%
and 
\begin{eqnarray*}
\varphi _{y}^{e}(t)-\varphi _{y}(t) &=&e\int_{0}^{1}\left\{ \varphi
_{xy}^{e,\alpha }(t)\left( x_{1}^{\varepsilon }(t)+x_{2}^{\varepsilon
}(t)\right) +\varphi _{yy}^{e,\alpha }(t)\mathbb{E}\left( x_{1}^{\varepsilon
}(t)+x_{2}^{\varepsilon }(t)\right) \right\} d\alpha \\
&&e\int_{0}^{1}\left\{ \varphi _{xy}^{e,\alpha }(t)\lambda ^{\varepsilon
}(t)+\varphi _{yy}^{e,\alpha }(t)\mathbb{E}\left( \lambda ^{\varepsilon
}(t)\right) \right\} d\alpha +\delta \varphi _{y}(t)\mathbf{I}_{\mathcal{E}%
_{\varepsilon }}(t).
\end{eqnarray*}%
Next we introduce the following notations:%
\begin{equation*}
\left\{ 
\begin{array}{l}
Z_{\varphi }^{1,\varepsilon }(t)=\dint_{0}^{1}\dint_{0}^{1}e\left\{ \varphi
_{xx}^{e,\alpha }(t)\lambda ^{\varepsilon }(t)\left( x_{1}^{\varepsilon
}(t)+x_{2}^{\varepsilon }(t)\right) \right. \smallskip \smallskip \\ 
\text{ \ \ \ \ \ \ \ \ \ }+\varphi _{xy}^{e,\alpha }(t)\left(
x_{1}^{\varepsilon }(t)+x_{2}^{\varepsilon }(t)\right) \mathbb{E}\left(
\lambda ^{\varepsilon }(t)\right) \left. +\lambda ^{\varepsilon }(t)\mathbb{E%
}\left( x_{1}^{\varepsilon }(t)+x_{2}^{\varepsilon }(t)\right) \right\}
d\alpha de\smallskip \smallskip \\ 
\text{ \ \ \ \ \ \ \ \ \ }+\dint_{0}^{1}\dint_{0}^{1}e\left\{ \varphi
_{yy}^{e,\alpha }(t)\mathbb{E}\left( \lambda ^{\varepsilon }(t)\right) 
\mathbb{E}\left( x_{1}^{\varepsilon }(t)+x_{2}^{\varepsilon }(t)\right)
\right\} d\alpha de\smallskip \smallskip \\ 
Z_{\varphi }^{2,\varepsilon }(t)=\dint_{0}^{1}\dint_{0}^{1}e\left\{ \varphi
_{xx}^{e,\alpha }(t)\left[ \left( x_{1}^{\varepsilon }(t)+x_{2}^{\varepsilon
}(t)\right) ^{2}-\left( x_{1}^{\varepsilon }(t)\right) ^{2}\right] \right.
\smallskip \smallskip \\ 
\text{ \ \ \ \ \ \ \ \ \ }\left. +2\varphi _{xy}^{e,\alpha }(t)\left(
x_{1}^{\varepsilon }(t)+x_{2}^{\varepsilon }(t)\right) \mathbb{E}\left(
x_{1}^{\varepsilon }(t)+x_{2}^{\varepsilon }(t)\right) \right\} d\alpha
de\smallskip \smallskip \\ 
\text{ \ \ \ \ \ \ \ \ \ }+\dint_{0}^{1}\dint_{0}^{1}e\left\{ \varphi
_{yy}^{e,\alpha }(t)\left( \mathbb{E}\left( x_{1}^{\varepsilon
}(t)+x_{2}^{\varepsilon }(t)\right) \right) ^{2}\right\} d\alpha
de\smallskip \smallskip \\ 
Z_{\varphi }^{3,\varepsilon }(t)=\dint_{0}^{1}\dint_{0}^{1}e\left\{ \varphi
_{xx}^{e,\alpha }(t)-\varphi _{xx}^{e,\alpha }(t)\left( x_{1}^{\varepsilon
}(t)\right) ^{2}\right\} d\alpha de\smallskip \smallskip \\ 
Z_{\varphi }^{4,\varepsilon }(t)=\left[ \delta \varphi
_{x}(t)x_{2}^{\varepsilon }(t)+\delta \varphi _{y}(t)\mathbb{E}\left(
x_{1}^{\varepsilon }(t)+x_{2}^{\varepsilon }(t)\right) \right] \mathbf{I}_{%
\mathcal{E}_{\varepsilon }}(t).%
\end{array}%
\right.
\end{equation*}%
From (\ref{3.42}) we get 
\begin{eqnarray*}
\Pi _{\varphi }^{\varepsilon }\left( t\right) &=&Z_{\varphi }^{1,\varepsilon
}(t)+Z_{\varphi }^{2,\varepsilon }(t)+Z_{\varphi }^{3,\varepsilon
}(t)+Z_{\varphi }^{4,\varepsilon }(t) \\
&&+\int_{0}^{1}\left\{ \varphi _{x}^{e}(t)\lambda ^{\varepsilon }(t)+\varphi
_{y}^{e}(t)\mathbb{E}\left( \lambda ^{\varepsilon }(t)\right) \right\} de,
\end{eqnarray*}%
applying (\ref{3.40}) together with estimates (\ref{3.15}) and (\ref{3.17})
we get $k\geq 1.$%
\begin{equation}
\mathbb{E}\left[ \sup_{t\in \left[ sT\right] }\left\vert \lambda
^{\varepsilon }(t)\right\vert ^{2k}\right] \leq C_{k}\varepsilon ^{2k}\text{.%
}  \tag{3.44}  \label{3.44}
\end{equation}%
Combining estimates (\ref{3.44}), (\ref{3.13}) and (\ref{3.15}) we get 
\begin{equation}
\mathbb{E}\left[ \sup_{t\in \left[ sT\right] }\left\vert Z_{\varphi
}^{1,\varepsilon }(t)\right\vert ^{2k}\right] \leq C_{k}\varepsilon ^{3k}%
\text{.}  \tag{3.45}  \label{3.45}
\end{equation}%
Similar arguments developed above with the helps of estimates (\ref{3.13}), (%
\ref{3.14}) and (\ref{3.15}) we can prove%
\begin{equation}
\mathbb{E}\left[ \sup_{t\in \left[ sT\right] }\left\vert Z_{\varphi
}^{2,\varepsilon }(t)\right\vert ^{2k}\right] \leq C_{k}\varepsilon
^{2k}\rho _{1,k}(\varepsilon )\text{.}  \tag{3.46}  \label{3.46}
\end{equation}%
where $\rho _{1,k}(\varepsilon )=\left( \varepsilon ^{k}+\varepsilon
^{2k}+\varepsilon ^{k}\rho ^{k}(\varepsilon )+\rho ^{k}(\varepsilon )\right)
\rightarrow 0$ as $\varepsilon \rightarrow 0.$ From \textit{Lebesgue's
bounded convergence theorem} it holds that%
\begin{equation}
\mathbb{E}\left[ \left( \int_{s}^{T}\left\vert Z_{\varphi }^{3,\varepsilon
}(t)\right\vert ^{2}dt\right) ^{k}\right] \leq C_{k}\varepsilon ^{2k}\left[ 
\mathbb{E}\left( \int_{s}^{T}\int_{0}^{1}\int_{0}^{1}\left\vert \varphi
_{xx}^{e,\alpha }(t)-\varphi _{xx}(t)\right\vert ^{4k}d\alpha dedt\right) %
\right] ^{\frac{1}{2}},  \tag{3.47}  \label{3.47}
\end{equation}%
here, if we denote $\rho _{2,k}(\varepsilon )=\left[ \mathbb{E}\left(
\int_{s}^{T}\int_{0}^{1}\int_{0}^{1}\left\vert \varphi _{xx}^{e,\alpha
}(t)-\varphi _{xx}(t)\right\vert ^{4k}d\alpha dedt\right) \right] ^{\frac{1}{%
2}}$ then $\lim_{\varepsilon \rightarrow 0}\rho _{2,k}(\varepsilon )=0.$
Also,%
\begin{equation}
\mathbb{E}\left[ \left( \int_{s}^{T}\left\vert Z_{\varphi }^{4,\varepsilon
}(t)\right\vert ^{2}\right) ^{k}\right] dt\leq C_{k}\varepsilon ^{2k}\rho
_{3,k}(\varepsilon ),  \tag{3.48}  \label{3.48}
\end{equation}%
where $\rho _{3,k}(\varepsilon )=\left( \varepsilon ^{k}+\rho
^{k}(\varepsilon )\right) \rightarrow 0,$ as $\varepsilon \rightarrow 0.$
Combining estimates (\ref{3.45})$\sim $(\ref{3.48}) we deduce%
\begin{equation}
\mathbb{E}\left[ \left( \int_{s}^{t}\left\vert \Pi _{\varphi }^{\varepsilon
}\left( r\right) \right\vert ^{2}dr\right) ^{k}\right] \leq C_{k}\varepsilon
^{2k}\rho _{k}(\varepsilon )+C_{k}\left[ \int_{s}^{t}\mathbb{E}\left(
\left\vert \lambda ^{\varepsilon }(r)\right\vert ^{2k}\right) dr\right] , 
\tag{3.49}  \label{3.49}
\end{equation}%
where $\rho _{k}(\varepsilon )=\left( \varepsilon ^{k}+\rho
_{1,k}(\varepsilon )+\rho _{2,k}(\varepsilon )+\rho _{3,k}(\varepsilon
)\right) \rightarrow 0,$ as $\varepsilon \rightarrow 0.\medskip $

\noindent Now, let us turn to estimate the jump terms $\Lambda
_{g}^{\varepsilon }\left( t,\theta \right) .$

\noindent \textit{Estimates of }$\Lambda _{g}^{\varepsilon }\left( t,\theta
\right) :$ We have for $t\in \left[ s,T\right] ,$%
\begin{equation*}
g(t,x^{\varepsilon }(t),u^{\varepsilon }(t),\theta )-\varphi (t,x^{\ast
}(t),u^{\ast }(t),\theta )=\int_{0}^{1}\left( g_{x}^{e}(t,\theta )\left(
x^{\varepsilon }(t)-x^{\ast }(t)\right) \right) de,
\end{equation*}%
where, for the subscript $\varkappa $ which indicates the first and the
second order derivatives of $g$, respectively, with respect to $\varkappa
=x,xx$, and for real\ $\hslash \in \left[ 0,1\right] :$%
\begin{equation*}
g_{\varkappa }^{\hslash }(t,\theta )=g_{\varkappa }(t,x^{\ast }(t)+\hslash
\left( x^{\varepsilon }(t)-x^{\ast }(t)\right) ,u^{\ast }(t),\theta ).
\end{equation*}%
Moreover,%
\begin{eqnarray*}
&&g(t,x^{\varepsilon }(t),u^{\varepsilon }(t),\theta )-g(t,x^{\ast
}(t),u^{\ast }(t),\theta )-g_{x}\left( t,\theta \right) \left(
x_{1}^{\varepsilon }(t)+x_{2}^{\varepsilon }(t)\right) \\
&=&\int_{0}^{1}\left\{ g_{x}^{e}(t,\theta )\lambda ^{\varepsilon }(t)+\left(
g_{x}^{e}(t,\theta )-g_{x}(t,\theta )\right) \left( x_{1}^{\varepsilon
}(t)+x_{2}^{\varepsilon }(t)\right) \right\} de.
\end{eqnarray*}%
By similar arguments we get%
\begin{eqnarray*}
g_{x}^{e}(t,\theta )-g_{x}(t,\theta ) &=&e\int_{0}^{1}\left\{
g_{xx}^{e,\alpha }(t,\theta )\left( x^{\varepsilon }(t)-x^{\ast }(t)\right) +%
\mathbb{E}\left( x^{\varepsilon }(t)-x^{\ast }(t)\right) \right\} d\alpha \\
&&+\delta g_{x}(t,\theta )\mathbf{I}_{\mathcal{E}_{\varepsilon }}(t) \\
&=&e\int_{0}^{1}\left\{ g_{xx}^{e,\alpha }(t,\theta )\lambda ^{\varepsilon
}(t)\right\} d\alpha \\
&&+e\int_{0}^{1}\left\{ g_{xx}^{e,\alpha }(t,\theta )\left(
x_{1}^{\varepsilon }(t)+x_{2}^{\varepsilon }(t)\right) \right\} d\alpha
+\delta g_{x}(t,\theta )\mathbf{I}_{\mathcal{E}_{\varepsilon }}(t).
\end{eqnarray*}%
Next we introduce the following notations:%
\begin{equation*}
\left\{ 
\begin{array}{l}
Z_{g}^{1,\varepsilon }(t,\theta )=\dint_{0}^{1}\dint_{0}^{1}e\left\{
g_{xx}^{e,\alpha }(t,\theta )\lambda ^{\varepsilon }(t)\left(
x_{1}^{\varepsilon }(t)+x_{2}^{\varepsilon }(t)\right) \right. \smallskip
\smallskip \\ 
\text{ \ \ \ \ \ \ \ \ \ \ \ \ \ \ \ \ }\left. +\lambda ^{\varepsilon }(t)%
\mathbb{E}\left( x_{1}^{\varepsilon }(t)+x_{2}^{\varepsilon }(t)\right)
\right\} d\alpha de\smallskip \smallskip \\ 
Z_{g}^{2,\varepsilon }(t,\theta )=\dint_{0}^{1}\dint_{0}^{1}e\left\{
g_{xx}^{e,\alpha }(t,\theta )\left[ \left( x_{1}^{\varepsilon
}(t)+x_{2}^{\varepsilon }(t)\right) ^{2}-\left( x_{1}^{\varepsilon
}(t)\right) ^{2}\right] \right. \smallskip \smallskip \\ 
Z_{g}^{3,\varepsilon }(t,\theta )=\dint_{0}^{1}\dint_{0}^{1}e\left\{
g_{xx}^{e,\alpha }(t,\theta )-g_{xx}^{e,\alpha }(t,\theta )\left(
x_{1}^{\varepsilon }(t)\right) ^{2}\right\} d\alpha de\smallskip \smallskip
\\ 
Z_{g}^{4,\varepsilon }(t,\theta )=\delta g_{x}(t,\theta )x_{2}^{\varepsilon
}(t)\mathbf{I}_{\mathcal{E}_{\varepsilon }}(t).%
\end{array}%
\right.
\end{equation*}%
From (\ref{3.43}) we get 
\begin{eqnarray*}
\Lambda _{g}^{\varepsilon }\left( t,\theta \right) &=&Z_{\varphi
}^{1,\varepsilon }(t,\theta )+Z_{\varphi }^{2,\varepsilon }(t,\theta
)+Z_{\varphi }^{3,\varepsilon }(t,\theta )+Z_{\varphi }^{4,\varepsilon
}(t,\theta ) \\
&&+\int_{0}^{1}g_{x}^{e}(t,\theta )\lambda ^{\varepsilon }(t)de.
\end{eqnarray*}%
By applying similar arguments developed in estimate $\Pi _{\varphi
}^{\varepsilon }\left( t\right) $ we can get%
\begin{equation}
\mathbb{E}\left[ \sup_{t\in \left[ sT\right] }\left\vert \Lambda
_{g}^{\varepsilon }\left( t,\theta \right) \right\vert ^{2k}\right] \leq
C_{k}\varepsilon ^{2k}\rho _{k}(\varepsilon ),  \tag{3.50}  \label{3.50}
\end{equation}%
where $\rho _{k}(\varepsilon )\rightarrow 0,$ as $\varepsilon \rightarrow 0.$

\noindent Finally by combining (\ref{3.49}), (\ref{3.50}) and (\ref{3.41})
with the help of Propositions A1, and Gronwall's Lemma, we conclude 
\begin{equation}
\mathbb{E}\left[ \sup_{t\in \left[ sT\right] }\left\vert \lambda
^{\varepsilon }(t)\right\vert ^{2k}\right] \leq C_{k,\mu (\Theta
)}\varepsilon ^{2k}\rho _{k}(\varepsilon ).  \tag{3.51}  \label{3.51}
\end{equation}%
This completes the proof of estimate (\ref{3.18}).

\noindent Noting that estimates (\ref{3.13}), (\ref{3.15}), (\ref{3.16}) and
(\ref{3.17}),\ follows from standard arguments.

\noindent Now by applying estimates (\ref{3.51}), (\ref{3.49}) the following
estimates hold.$\smallskip $

\noindent \textbf{Corollary 3.1}. We have for $\varphi =f,\sigma ,\ell $ 
\begin{equation}
\mathbb{E}\left[ \left( \int_{s}^{T}\left\vert \Pi _{\varphi }^{\varepsilon
}\left( r\right) \right\vert ^{2}dr\right) ^{k}\right] \leq C_{k}\varepsilon
^{2k}\rho _{k}(\varepsilon ),  \tag{3.52}  \label{3.52}
\end{equation}%
\begin{equation}
\mathbb{E}\left[ \left\vert \Pi _{h}^{\varepsilon }\left( T\right)
\right\vert \right] \leq C_{k}\varepsilon \rho (\varepsilon ),  \tag{3.53}
\label{3.53}
\end{equation}%
where 
\begin{eqnarray*}
\Pi _{h}^{\varepsilon }\left( T\right) &=&h\left( x^{\varepsilon }(T),%
\mathbb{E}(x^{\varepsilon }(T))\right) -h(x^{\ast }(T),\mathbb{E}(x^{\ast
}(T)))-h_{x}\left( T\right) \left( x_{1}^{\varepsilon
}(T)+x_{2}^{\varepsilon }(T)\right) \\
&&-\left\{ h_{y}\left( T\right) \mathbb{E}\left( x_{1}^{\varepsilon
}(T)+x_{2}^{\varepsilon }(T)\right) +\mathcal{L}_{T}\left(
h,x_{1}^{\varepsilon }\right) \right\} ,
\end{eqnarray*}
and $\rho _{k}(\varepsilon ),$ $\rho (\varepsilon )$ tends to $0$\ as $%
\varepsilon \rightarrow 0.\smallskip $

\noindent \textbf{Lemma 3.4. }We have%
\begin{equation}
\begin{array}{l}
\mathbb{E}\left[ h_{xx}\left( x^{\ast }(T),\mathbb{E}(x^{\ast }(T))\right)
x_{1}^{\varepsilon }(T)^{2}\right] \smallskip \smallskip \\ 
=-\mathbb{E}\dint_{s}^{T}\left\{ -H_{xx}(t)\left( x_{1}^{\varepsilon
}(t)\right) ^{2}\right. +Q^{\ast }(t)\sigma _{y}^{2}(t)\left( \mathbb{E}%
\left( x_{1}^{\varepsilon }(t)\right) \right) ^{2}+Q^{\ast }(t)(\left(
\delta \sigma (t)\right) ^{2}\smallskip \smallskip \\ 
\text{ \ \ \ \ }+\dint_{\Theta }\left( \delta g(t,\theta \right) )^{2}\mu
(d\theta ))\mathbf{I}_{\mathcal{E}_{\varepsilon }}(t)+\Gamma _{t}^{\ast
}(\theta )\dint_{\Theta }\left( \delta g(t,\theta \right) )^{2}\mathbf{I}_{%
\mathcal{E}_{\varepsilon }}(t)\mu (d\theta )\smallskip \smallskip \\ 
\text{ \ \ \ \ }+2\left( \mathbb{E}\left( x_{1}^{\varepsilon }(t)\right)
\right) x_{1}^{\varepsilon }(t)\left[ Q^{\ast }(t)f_{y}(t)+Q^{\ast
}(t)\sigma _{x}(t)\sigma _{y}(t)+R^{\ast }(t)\sigma _{y}(t)\right]
\smallskip \smallskip \\ 
\text{ \ \ \ \ }+\left. Q^{\ast }(t)\dint_{\Theta }\left( g_{x}\left(
t,\theta \right) \right) ^{2}\left( x_{1}^{\varepsilon }(t)\right) ^{2}\mu
(d\theta )\right\} dt.%
\end{array}
\tag{3.54}  \label{3.54}
\end{equation}

\noindent \textbf{Proof.} By using integration by parts formula for jumps
processes to $Q^{\ast }(t)\left( x_{1}^{\varepsilon }(t)\right) ^{2}$ (see
Lemma A1) and taking expectation, we get from (\ref{3.13}) and (\ref{3.14})%
\begin{equation}
\begin{array}{l}
\mathbb{E(}Q^{\ast }(T)x_{1}^{\varepsilon }(T)^{2})=\mathbb{E}%
\dint_{s}^{T}Q^{\ast }(t)d(\left( x_{1}^{\varepsilon }(t)\right) ^{2})+%
\mathbb{E}\dint_{s}^{T}\left( x_{1}^{\varepsilon }(t)\right) ^{2}dQ^{\ast
}(t)\smallskip \smallskip \\ 
\text{ \ \ \ }+\mathbb{E}\dint_{s}^{T}R^{\ast }(t)2x_{1}^{\varepsilon }(t)%
\left[ \sigma _{x}(t)x_{1}^{\varepsilon }(t)+\sigma _{y}(t)\mathbb{E}\left(
x_{1}^{\varepsilon }(t)\right) +\delta \sigma (t)\mathbf{I}_{\mathcal{E}%
_{\varepsilon }}(t)\right] dt\smallskip \smallskip \\ 
\text{ \ \ \ }+\mathbb{E}\dint_{s}^{T}\dint_{\Theta }\Gamma _{t}^{\ast
}(\theta )2x_{1}^{\varepsilon }(t)\left[ g_{x}\left( t,\theta \right)
x_{1}^{\varepsilon }(t)+\delta g(t,\theta )\mathbf{I}_{\mathcal{E}%
_{\varepsilon }}(t)\right] \mu \left( d\theta \right) dt\smallskip \smallskip
\\ 
\text{ \ \ \ \ }=\mathcal{J}_{1}^{\varepsilon }\mathbb{+}\mathcal{J}_{2}^{%
\mathbb{\varepsilon }}\mathbb{+}\mathcal{J}_{3}^{\varepsilon }\mathbb{+}%
\mathcal{J}_{4}^{\varepsilon }.%
\end{array}
\tag{3.55}  \label{3.55}
\end{equation}%
By using It\^{o} formula to jump process $\left( x_{1}^{\varepsilon
}(t)\right) ^{2}$ (see Situ \cite{situ}) we have%
\begin{equation}
\begin{array}{l}
\mathcal{J}_{1}^{\varepsilon }=\mathbb{E}\dint_{s}^{T}Q^{\ast }(t)d(\left(
x_{1}^{\varepsilon }(t)\right) ^{2})\smallskip \smallskip \\ 
\text{ \ \ \ }=\mathbb{E}\dint_{s}^{T}Q^{\ast }(t)\left\{
2x_{1}^{\varepsilon }(t)\left[ f_{x}(t)x_{1}^{\varepsilon }(t)+f_{y}(t)%
\mathbb{E}\left( x_{1}^{\varepsilon }(t)\right) +\delta f(t)\mathbf{I}_{%
\mathcal{E}_{\varepsilon }}(t)\right] \right. \smallskip \smallskip \\ 
\text{ \ \ \ }+\left\{ \sigma _{x}(t)x_{1}^{\varepsilon }(t)+\sigma _{y}(t)%
\mathbb{E}\left( x_{1}^{\varepsilon }(t)\right) +\delta \sigma (t)\mathbf{I}%
_{\mathcal{E}_{\varepsilon }}(t)\right\} ^{2}\smallskip \smallskip \\ 
\text{ \ \ \ }+\left. \dint_{\Theta }\left\{ g_{x}\left( t,\theta \right)
x_{1}^{\varepsilon }(t)+\delta g(t,\theta )\mathbf{I}_{\mathcal{E}%
_{\varepsilon }}(t)\right\} ^{2}\mu (d\theta )\right\} dt.\smallskip
\smallskip%
\end{array}
\tag{3.56}  \label{3.56}
\end{equation}%
Applying (\ref{2.6}) we can get

\begin{equation}
\begin{array}{l}
\mathcal{J}_{2}^{\varepsilon }=\mathbb{E}\dint_{s}^{T}\left(
x_{1}^{\varepsilon }(t)\right) ^{2}dQ^{\ast }(t)\smallskip \smallskip \\ 
\text{ \ \ \ }=-\mathbb{E}\dint_{s}^{T}\left( x_{1}^{\varepsilon }(t)\right)
^{2}\left\{ 2f_{x}\left( t\right) Q^{\ast }(t)+\sigma _{x}^{2}\left(
t\right) Q^{\ast }(t)+2\sigma _{x}\left( t\right) R^{\ast }(t)\right.
\smallskip \smallskip \\ 
\text{ \ \ \ }+\dint_{\Theta }\left( g_{x}\left( t,\theta \right) \right)
^{2}\left( \Gamma _{t}^{\ast }(\theta )+Q^{\ast }(t)\right) \mu (d\theta
)+2\dint_{\Theta }\Gamma _{t}^{\ast }(\theta )g_{x}\left( t,\theta \right)
\mu (d\theta )\smallskip \smallskip \\ 
\text{ \ \ \ }+\left. H_{xx}(t))\right\} dt.%
\end{array}
\tag{3.57}  \label{3.57}
\end{equation}%
A simple computations shows that

\begin{equation}
\begin{array}{l}
\mathcal{J}_{3}^{\varepsilon }=\mathbb{E}\dint_{s}^{T}R^{\ast
}(t)2x_{1}^{\varepsilon }(t)\left[ \sigma _{x}(t)x_{1}^{\varepsilon
}(t)+\sigma _{y}(t)\mathbb{E}\left( x_{1}^{\varepsilon }(t)\right) +\delta
\sigma (t)\mathbf{I}_{\mathcal{E}_{\varepsilon }}(t)\right] dt\smallskip
\smallskip \\ 
\text{ \ \ }=2\mathbb{E}\dint_{s}^{T}\left\{ R^{\ast }(t)\sigma
_{x}(t)\left( x_{1}^{\varepsilon }(t)\right) ^{2}\right. \smallskip
\smallskip \\ 
\text{ \ \ }+R^{\ast }(t)\sigma _{y}(t)\mathbb{E}\left( x_{1}^{\varepsilon
}(t)\right) x_{1}^{\varepsilon }(t)dt\smallskip \smallskip \\ 
\text{ \ \ }\left. +R^{\ast }(t)\delta \sigma (t)x_{1}^{\varepsilon }(t)%
\mathbf{I}_{\mathcal{E}_{\varepsilon }}(t)\right\} dt,%
\end{array}
\tag{3.58}  \label{3.58}
\end{equation}%
and 
\begin{equation}
\begin{array}{l}
\mathcal{J}_{4}^{\varepsilon }=2\mathbb{E}\dint_{s}^{T}\dint_{\Theta }\Gamma
_{t}^{\ast }(\theta )x_{1}^{\varepsilon }(t)\left[ g_{x}\left( t,\theta
\right) x_{1}^{\varepsilon }(t)+\delta g(t,\theta )\mathbf{I}_{\mathcal{E}%
_{\varepsilon }}(t)\right] \mu \left( d\theta \right) dt\smallskip \smallskip
\\ 
\text{ \ \ }=2\mathbb{E}\dint_{s}^{T}\dint_{\Theta }\Gamma _{t}^{\ast
}(\theta )g_{x}\left( t,\theta \right) \left( x_{1}^{\varepsilon }(t)\right)
^{2}\mu \left( d\theta \right) dt\smallskip \smallskip \\ 
\text{ \ \ }+2\mathbb{E}\dint_{s}^{T}\dint_{\Theta }\Gamma _{t}^{\ast
}(\theta )\delta g(t,\theta )x_{1}^{\varepsilon }(t)\mathbf{I}_{\mathcal{E}%
_{\varepsilon }}(t)\mu \left( d\theta \right) dt.%
\end{array}
\tag{3.59}  \label{3.59}
\end{equation}%
Thus, by combining (\ref{3.56})$\sim $(\ref{3.59}) together with (\ref{3.55}%
) it follows that

\begin{equation*}
\begin{array}{l}
\mathbb{E(}Q^{\ast }(T)\left( x_{1}^{\varepsilon }(T)\right) ^{2})\smallskip
\smallskip \\ 
\text{ \ \ \ }=\mathbb{E}\dint_{s}^{T}\left\{ -H_{xx}(t)\left(
x_{1}^{\varepsilon }(t)\right) ^{2}\right. \smallskip \smallskip \\ 
\text{ \ \ }+Q^{\ast }(t)\sigma _{y}^{2}(t)\left( \mathbb{E}\left(
x_{1}^{\varepsilon }(t)\right) \right) ^{2}+Q^{\ast }(t)\left( \left( \delta
\sigma (t)\right) ^{2}+\dint_{\Theta }\left( \delta g(t,\theta \right)
)^{2}\mu (d\theta )\right) \mathbf{I}_{\mathcal{E}_{\varepsilon
}}(t)\smallskip \smallskip \\ 
\text{ \ \ }+\dint_{\Theta }\Gamma _{t}^{\ast }(\theta )\left( \delta
g(t,\theta \right) )^{2}\mathbf{I}_{\mathcal{E}_{\varepsilon }}(t)\mu
(d\theta )\smallskip \smallskip \\ 
\text{ \ \ }+2\left( \mathbb{E}\left( x_{1}^{\varepsilon }(t)\right) \right)
x_{1}^{\varepsilon }(t)\left[ Q^{\ast }(t)f_{y}(t)+Q^{\ast }(t)\sigma
_{x}(t)\sigma _{y}(t)+R^{\ast }(t)\sigma _{y}(t)\right] \smallskip \smallskip
\\ 
\text{ \ \ }+\left. Q^{\ast }(t)\dint_{\Theta }\left( g_{x}\left( t,\theta
\right) \right) ^{2}\left( x_{1}^{\varepsilon }(t)\right) ^{2}\mu (d\theta
)\right\} dt.%
\end{array}%
\end{equation*}%
Finally, since $Q^{\ast }(T)=-h_{xx}\left( x^{\ast }(T),\mathbb{E}(x^{\ast
}(T))\right) ,$ this completes the proof of\textbf{\ }Lemma 3.4.$\smallskip $

\noindent The following Lemma gives estimates related to the adjoint
processes $(\Psi ^{\ast }\left( \cdot \right) ,K^{\ast }\left( \cdot \right)
,$ $\gamma ^{\ast }\left( \cdot \right) )$\ and $\left( Q^{\ast }\left(
\cdot \right) ,R^{\ast }\left( \cdot \right) ,\Gamma ^{\ast }\left( \cdot
\right) \right) $\ given by (\ref{2.5})$,$ (\ref{2.6}) respectively.

\noindent \textbf{Lemma 3.5. }We have

\begin{equation}
\begin{array}{c}
\mathbb{E}\left\{ \dint_{s}^{T}\left\vert \left[ \Psi ^{\ast }\left(
t\right) \delta f_{x}(t)+K^{\ast }\left( t\right) \delta \sigma
_{x}(t)+\dint_{\Theta }\gamma _{t}^{\ast }\left( \theta \right) \delta
g_{x}(t,\theta )\mu (d\theta )\right] x_{1}^{\varepsilon }(t)\mathbf{I}_{%
\mathcal{E}_{\varepsilon }}(t)\right\vert dt\right\} \smallskip \smallskip
\\ 
\leq C\varepsilon \rho (\varepsilon ),%
\end{array}
\tag{3.60}  \label{3.60}
\end{equation}%
\begin{equation}
\begin{array}{c}
\mathbb{E}\left\{ \dint_{s}^{T}\left\vert \left[ Q^{\ast }\left( t\right)
f_{y}(t)+Q^{\ast }\left( t\right) \sigma _{x}(t)\sigma _{y}(t)+R^{\ast
}(t)\sigma _{y}(t)\right] x_{1}^{\varepsilon }(t)\mathbb{E}\left(
x_{1}^{\varepsilon }(t)\right) \right\vert dt\right\} \smallskip \smallskip
\\ 
\leq C\varepsilon \rho (\varepsilon ),%
\end{array}
\tag{3.61}  \label{3.61}
\end{equation}%
and%
\begin{equation}
\mathbb{E}\left\{ \dint_{s}^{T}\left\vert Q^{\ast }\left( t\right) \left(
\sigma _{y}(t)\right) ^{2}\left( \mathbb{E}\left( x_{1}^{\varepsilon
}(t)\right) \right) ^{2}\right\vert dt\right\} \leq C\varepsilon \rho
(\varepsilon ),  \tag{3.62}  \label{3.62}
\end{equation}%
\begin{equation}
\mathbb{E}\left\{ \int_{s}^{T}\left\vert \int_{\Theta }Q^{\ast }(t)\left(
g_{x}\left( t,\theta \right) \right) ^{2}\left( x_{1}^{\varepsilon
}(t)\right) ^{2}\mu (d\theta )\right\vert dt\right\} \leq C\varepsilon \rho
(\varepsilon ),  \tag{3.63}  \label{3.63}
\end{equation}%
where $\rho (\varepsilon )\rightarrow 0$ as $\varepsilon \rightarrow 0.$

\noindent \textbf{Proof. }

\noindent \textit{Estimates of (\ref{3.60}):} First we have 
\begin{equation}
\begin{array}{l}
\mathbb{E}\left\{ \dint_{s}^{T}\left\vert \left[ \Psi ^{\ast }\left(
t\right) \delta f_{x}(t)+K^{\ast }\left( t\right) \delta \sigma
_{x}(t)+\dint_{\Theta }\gamma _{t}^{\ast }\left( \theta \right) \delta
g_{x}(t,\theta )\mu (d\theta )\right] x_{1}^{\varepsilon }(t)\mathbf{I}_{%
\mathcal{E}_{\varepsilon }}(t)\right\vert dt\right\} \smallskip \smallskip
\\ 
\text{ \ \ \ \ \ \ \ \ }\leq \mathbb{E}\left[ \dint_{s}^{T}\left\vert \Psi
^{\ast }\left( t\right) \delta f_{x}(t)x_{1}^{\varepsilon }(t)\mathbf{I}_{%
\mathcal{E}_{\varepsilon }}(t)\right\vert dt\right] +\mathbb{E}\left[
\dint_{s}^{T}\left\vert K^{\ast }\left( t\right) \delta \sigma
_{x}(t)x_{1}^{\varepsilon }(t)\mathbf{I}_{\mathcal{E}_{\varepsilon
}}(t)\right\vert dt\right] \smallskip \smallskip \\ 
\text{ \ \ \ \ \ \ \ \ }+\mathbb{E}\left[ \dint_{s}^{T}\left\vert
\dint_{\Theta }\gamma _{t}^{\ast }\left( \theta \right) \delta
g_{x}(t,\theta )\mu (d\theta )x_{1}^{\varepsilon }(t)\mathbf{I}_{\mathcal{E}%
_{\varepsilon }}(t)\right\vert dt\right] \smallskip \smallskip \\ 
\text{ \ \ \ \ \ \ \ \ }=\mathcal{I}_{1}^{\varepsilon }+\mathcal{I}%
_{2}^{\varepsilon }+\mathcal{I}_{3}^{\varepsilon }.%
\end{array}
\tag{3.64}  \label{3.64}
\end{equation}%
Using (\ref{2.7}) and estimates (\ref{3.13} with $k=1$), then from
Cauchy-Schwarz inequality we get%
\begin{equation}
\begin{array}{l}
\mathcal{I}_{2}^{\varepsilon }=\mathbb{E}\left[ \dint_{s}^{T}\left\vert
K^{\ast }\left( t\right) \delta \sigma _{x}(t)x_{1}^{\varepsilon }(t)\mathbf{%
I}_{\mathcal{E}_{\varepsilon }}(t)\right\vert dt\right] \smallskip \smallskip
\\ 
\text{ \ \ \ \ \ }\leq C\left[ \mathbb{E}\left( \sup_{t\in \left[ s,T\right]
}\left\vert x_{1}^{\varepsilon }(t)\right\vert ^{2}\right) \right] ^{\frac{1%
}{2}}\left[ \mathbb{E}\left( \left( \dint_{s}^{T}\left\vert K^{\ast }\left(
t\right) \right\vert \mathbf{I}_{\mathcal{E}_{\varepsilon }}(t)dt\right)
^{2}\right) \right] ^{\frac{1}{2}}\smallskip \smallskip \\ 
\text{ \ \ \ \ \ }\leq C\mathbb{\varepsilon }^{\frac{1}{2}}\left[ \mathbb{E}%
\left( \left( \dint_{s}^{T}\left\vert K^{\ast }\left( t\right) \right\vert 
\mathbf{I}_{\mathcal{E}_{\varepsilon }}(t)dt\right) ^{2}\right) \right] ^{%
\frac{1}{2}}\smallskip \smallskip \\ 
\text{ \ \ \ \ \ }\leq C\mathbb{\varepsilon }^{\frac{1}{2}}\left[ \mathbb{E}%
\left( \dint_{s}^{T}\left\vert K^{\ast }\left( t\right) \right\vert ^{2}%
\mathbf{I}_{\mathcal{E}_{\varepsilon }}(t)dt\right) \right] ^{\frac{1}{2}}%
\mathbb{\varepsilon }^{\frac{1}{2}}\leq C\mathbb{\varepsilon \rho }%
_{2}\left( \varepsilon \right) ,%
\end{array}
\tag{3.65}  \label{3.65}
\end{equation}%
where, also from (\ref{2.7}) and\textit{\ Dominated Convergence Theorem }we
obtain 
\begin{equation*}
\mathbb{\rho }_{2}\left( \varepsilon \right) =\left[ \mathbb{E}\left(
\int_{s}^{T}\left\vert K^{\ast }\left( t\right) \right\vert ^{2}\mathbf{I}_{%
\mathcal{E}_{\varepsilon }}(t)dt\right) \right] ^{\frac{1}{2}}\rightarrow 0%
\text{ as }\varepsilon \rightarrow 0.
\end{equation*}

\noindent Similarly, we can prove estimate $\mathcal{I}_{1}^{\varepsilon }$
then we get%
\begin{equation}
\mathcal{I}_{1}^{\varepsilon }\leq C\mathbb{\varepsilon \rho }_{1}\left(
\varepsilon \right) .  \tag{3.66}  \label{3.66}
\end{equation}%
Let us turn to third term $\mathcal{I}_{3}^{\varepsilon }$. By using (\ref%
{2.7}) and estimate (\ref{3.13} with $k=1$) with the help of Cauchy-Schwarz
inequality we get

\begin{equation}
\begin{array}{l}
\mathcal{I}_{3}^{\varepsilon }=\mathbb{E}\left[ \dint_{s}^{T}\left\vert
\dint_{\Theta }\gamma _{t}^{\ast }\left( \theta \right) \delta
g_{x}(t,\theta )\mu (d\theta )x_{1}^{\varepsilon }(t)\mathbf{I}_{\mathcal{E}%
_{\varepsilon }}(t)\right\vert dt\right] \smallskip \smallskip \\ 
\ \ \ \ \leq C\left[ \mathbb{E}\left( \sup_{t\in \left[ s,T\right]
}\left\vert x_{1}^{\varepsilon }(t)\right\vert ^{2}\right) \right] ^{\frac{1%
}{2}}\left[ \mathbb{E}\left( \left( \dint_{s}^{T}\dint_{\Theta }\left\vert
\gamma _{t}^{\ast }\left( \theta \right) \right\vert \mathbf{I}_{\mathcal{E}%
_{\varepsilon }}(t)\mu (d\theta )dt\right) ^{2}\right) \right] ^{\frac{1}{2}%
}\smallskip \smallskip \\ 
\ \ \ \ \leq C\mathbb{\varepsilon }^{\frac{1}{2}}\left[ \mathbb{E}\left(
\dint_{s}^{T}\dint_{\Theta }\left\vert \gamma _{t}^{\ast }\left( \theta
\right) \right\vert ^{2}\mathbf{I}_{\mathcal{E}_{\varepsilon }}(t)\mu
(d\theta )dt\right) \right] ^{\frac{1}{2}}\smallskip \smallskip \\ 
\ \ \ \ \leq C\mu (\Theta )\mathbb{\varepsilon }^{\frac{1}{2}}\left[ \mathbb{%
E}\left( \dint_{s}^{T}\sup_{\theta \in \Theta }\left\vert \gamma _{t}^{\ast
}\left( \theta \right) \right\vert ^{2}\mathbf{I}_{\mathcal{E}_{\varepsilon
}}(t)dt\right) \right] ^{\frac{1}{2}}\mathbb{\varepsilon }^{\frac{1}{2}}\leq
C\mathbb{\varepsilon \rho }_{3}\left( \varepsilon \right) ,%
\end{array}
\tag{3.67}  \label{3.67}
\end{equation}%
Again, from (\ref{2.7}) and\textit{\ Dominated Convergence Theorem }we
obtain 
\begin{equation*}
\mathbb{\rho }_{3}\left( \varepsilon \right) =\left[ \mathbb{E}\left(
\int_{s}^{T}\sup_{\theta \in \Theta }\left\vert \gamma _{t}^{\ast }\left(
\theta \right) \right\vert ^{2}\mathbf{I}_{\mathcal{E}_{\varepsilon
}}(t)dt\right) \right] ^{\frac{1}{2}}\rightarrow 0\text{ as }\varepsilon
\rightarrow 0.
\end{equation*}%
Finally, we set $\mathbb{\rho }\left( \varepsilon \right) =\mathbb{\rho }%
_{1}\left( \varepsilon \right) +\mathbb{\rho }_{2}\left( \varepsilon \right)
+\mathbb{\rho }_{3}\left( \varepsilon \right) \rightarrow 0$ as $\varepsilon
\rightarrow 0$ then the desired result follows immediately from combining (%
\ref{3.63})$\sim $(\ref{3.67}). This completes the proof of (\ref{3.60}).

\noindent \textit{Estimates of (\ref{3.63}):} First we have from assumption
(H2) and by using (\ref{2.8}) and estimate (\ref{3.13} with $k=1$) with the
help of Cauchy-Schwarz inequality we get 
\begin{eqnarray*}
&&\mathbb{E}\left\{ \int_{s}^{T}\left\vert \int_{\Theta }Q^{\ast }(t)\left(
g_{x}\left( t,\theta \right) \right) ^{2}\left( x_{1}^{\varepsilon
}(t)\right) ^{2}\mu (d\theta )\right\vert dt\right\} \\
&\leq &C\mu (\Theta )\mathbb{E}\left\{ \int_{s}^{T}\left\vert Q^{\ast
}(t)\sup_{\theta \in \Theta }(g_{x}\left( t,\theta \right) ^{2}\left(
x_{1}^{\varepsilon }(t)\right) ^{2}\right\vert dt\right\} \\
&\leq &C\left[ \mathbb{E}\left( \sup \left\vert \left( x_{1}^{\varepsilon
}(t)\right) \right\vert ^{4}\right) \right] ^{\frac{1}{2}}\left[ \mathbb{E}%
\left( \left( \int_{s}^{T}\left\vert Q^{\ast }(t)\right\vert dt\right)
^{2}\right) \right] ^{\frac{1}{2}} \\
&\leq &C\varepsilon \rho (\varepsilon ),
\end{eqnarray*}%
where $\mathbb{\rho }\left( \varepsilon \right) \rightarrow 0$ as $%
\varepsilon \rightarrow 0$

\noindent Using similar arguments developed above for estimates (\ref{3.61})
and (\ref{3.62}) which completes the proof of Lemma 3.5.

\noindent It worth mentioning that by combining the duality relations (\ref%
{3.6}) and (\ref{3.7}) in Lemma 3.1 together with Lemma 3.5 we get

\begin{equation}
\begin{array}{l}
\mathbb{E}\left( \Psi (T)\left( x_{1}^{\varepsilon }(T)+x_{2}^{\varepsilon
}(T)\right) \right) =\mathbb{E}\dint_{s}^{T}\left( x_{1}^{\varepsilon
}(t)+x_{2}^{\varepsilon }(t)\right) \left[ \left( \ell _{x}(t)+\mathbb{E(}%
\ell _{y}(t)\right) \right] dt\smallskip \smallskip \\ 
\text{ \ }+\mathbb{E}\dint_{s}^{T}\left\{ \Psi (t)\delta f(t)+K(t)\delta
\sigma (t)+\dint_{\Theta }\gamma _{t}(\theta )\delta g(t,\theta )\mu \left(
d\theta \right) \right\} \mathbf{I}_{\mathcal{E}_{\varepsilon
}}(t)dt\smallskip \smallskip \\ 
\text{ \ }+\mathbb{E}\dint_{s}^{T}\left\{ \Psi (t)\mathcal{L}%
_{t}(f,x_{1}^{\varepsilon })+K(t)\mathcal{L}_{t}(\sigma ,x_{1}^{\varepsilon
})+\dint_{\Theta }\gamma _{t}(\theta )\mathcal{L}_{t,\theta
}(g,x_{1}^{\varepsilon })\mu \left( d\theta \right) \right\} dt\smallskip
\smallskip \\ 
\text{ \ }+\tau \left( \varepsilon \right) .%
\end{array}
\tag{3.68}  \label{3.68}
\end{equation}

\noindent \textbf{Proof of Theorem 3.1. }By applying (\ref{3.2}), (\ref{3.13}%
) and Corollary 3.1\textbf{\ }we get

\begin{eqnarray*}
0 &\leq &J^{^{s,\zeta }}\left( u^{\varepsilon }(\cdot )\right) -J^{^{s,\zeta
}}\left( u^{\ast }(\cdot )\right) \\
&=&\mathbb{E}\left[ h(x^{\varepsilon }(T),\mathbb{E}\left( x^{\varepsilon
}(T)\right) )-h(x^{\ast }(T),\mathbb{E}\left( x^{\ast }(T)\right) \right] \\
&&+\mathbb{E}\int_{s}^{T}\left[ \ell (t,x^{\varepsilon }(t),\mathbb{E}%
(x^{\varepsilon }(t)),u^{\varepsilon }(t))-\ell (t,x^{\ast }(t),\mathbb{E}%
(x^{\ast }(t)),u^{\ast }(t))\right] dt \\
&=&\mathbb{E}\left[ h_{x}(x^{\ast }(T),\mathbb{E}\left( x^{\ast }(T)\right)
)\left( x_{1}^{\varepsilon }(T)+x_{2}^{\varepsilon }(T)\right) \right] \\
&&+\mathbb{E}\left[ h_{y}(x^{\ast }(T),\mathbb{E}\left( x^{\ast }(T)\right)
\left( \mathbb{E}\left( x_{1}^{\varepsilon }(T)\right) +\mathbb{E}\left(
x_{2}^{\varepsilon }(T)\right) \right) \right] \\
&&+\mathbb{E}\int_{s}^{T}\left[ \ell _{x}\left( t\right) \left(
x_{1}^{\varepsilon }(t)+x_{2}^{\varepsilon }(t)\right) +\ell _{y}\left(
t\right) \left( \mathbb{E}\left( x_{1}^{\varepsilon }(t)\right) +\mathbb{E}%
\left( x_{2}^{\varepsilon }(t)\right) \right) \right] dt \\
&&+\mathbb{E}\int_{s}^{T}\left[ \delta \ell (t)\mathbf{I}_{\mathcal{E}%
_{\varepsilon }}(t)+\mathcal{L}_{t}\left( \ell ,x_{1}^{\varepsilon }\right) %
\right] dt+\mathbb{E}\left[ \mathcal{L}_{T}\left( h,x_{1}^{\varepsilon
}\right) \right] +\tau \left( \varepsilon \right) ,
\end{eqnarray*}%
then we get%
\begin{eqnarray*}
0 &\leq &J^{^{s,\zeta }}\left( u^{\varepsilon }(\cdot )\right) -J^{^{s,\zeta
}}\left( u^{\ast }(\cdot )\right) \\
&=&\mathbb{E}\int_{s}^{T}\left[ \delta \ell (t)\mathbf{I}_{\mathcal{E}%
_{\varepsilon }}(t)+\mathcal{L}_{t}\left( \ell ,x_{1}^{\varepsilon }\right) %
\right] dt+\mathbb{E}\left[ \mathcal{L}_{T}\left( h,x_{1}^{\varepsilon
}\right) \right] \\
&&+\mathbb{E}\int_{s}^{T}\left[ \ell _{x}\left( t\right) +\mathbb{E}\left(
\ell _{y}\left( t\right) \right) \right] \left( x_{1}^{\varepsilon
}(t)+x_{2}^{\varepsilon }(t)\right) dt \\
&&+\mathbb{E}\left\{ \left[ h_{x}(x^{\ast }(T),\mathbb{E}\left( x^{\ast
}(T)\right) +\mathbb{E}\left( h_{y}(x^{\ast }(T),\mathbb{E}\left( x^{\ast
}(T)\right) \right) \right] \left[ x_{1}^{\varepsilon
}(T)+x_{2}^{\varepsilon }(T)\right] \right\} +\tau \left( \varepsilon
\right) .
\end{eqnarray*}%
from (\ref{3.68}) and the fact that $\Psi ^{\ast }\left( T\right)
=-h_{x}(x^{\ast }(T),\mathbb{E}\left( x^{\ast }(T)\right) -\mathbb{E}\left(
h_{y}(x^{\ast }(T),\mathbb{E}\left( x^{\ast }(T)\right) \right) $ we obtain%
\begin{eqnarray*}
0 &\leq &J^{^{s,\zeta }}\left( u^{\varepsilon }(\cdot )\right) -J^{^{s,\zeta
}}\left( u^{\ast }(\cdot )\right) =\mathbb{E}\int_{s}^{T}\left[ \delta \ell
(t)\mathbf{I}_{\mathcal{E}_{\varepsilon }}(t)+\mathcal{L}_{t}\left( \ell
,x_{1}^{\varepsilon }\right) \right] dt+\mathbb{E}\left[ \mathcal{L}%
_{T}\left( h,x_{1}^{\varepsilon }\right) \right] \\
&&-\mathbb{E}\int_{s}^{T}\left\{ \Psi ^{\ast }\left( t\right) \delta
f(t)+K^{\ast }\left( t\right) \delta \sigma (t)+\int_{\Theta }\gamma
_{t}^{\ast }\left( \theta \right) \delta g(t,\theta )\mu \left( d\theta
\right) \right\} \mathbf{I}_{\mathcal{E}_{\varepsilon }}(t)dt \\
&&-\mathbb{E}\int_{s}^{T}\left\{ \Psi ^{\ast }\left( t\right) \mathcal{L}%
_{t}\left( f,x_{1}^{\varepsilon }\right) +K^{\ast }\left( t\right) \mathcal{L%
}_{t}\left( \sigma ,x_{1}^{\varepsilon }\right) \right. \\
&&+\left. \int_{\Theta }\gamma _{t}^{\ast }\left( \theta \right) \mathcal{L}%
_{t,\theta }\left( g,x_{1}^{\varepsilon }\right) \mu \left( d\theta \right)
\right\} dt+\tau \left( \varepsilon \right) .
\end{eqnarray*}%
Next by applying (\ref{2.9}) we deduce

\begin{equation}
\begin{array}{l}
0\leq J^{^{s,\zeta }}\left( u^{\varepsilon }(\cdot )\right) -J^{^{s,\zeta
}}\left( u^{\ast }(\cdot )\right) \smallskip \smallskip \\ 
\text{ \ \ \ \ \ }=-\mathbb{E}\dint_{s}^{T}\delta H(t)\mathbf{I}_{\mathcal{E}%
_{\varepsilon }}(t)dt\smallskip \smallskip \\ 
\text{ \ \ \ \ \ }+\frac{1}{2}\mathbb{E}\left[ h_{xx}\left( x^{\ast }(T),%
\mathbb{E}(x^{\ast }(T)\right) \left( x_{1}^{\varepsilon }(T)\right)
^{2}-\dint_{s}^{T}H_{xx}(t)\left( x_{1}^{\varepsilon }(t)\right) ^{2}dt%
\right] \smallskip \smallskip \\ 
\text{ \ \ \ \ \ }+\tau \left( \varepsilon \right) .%
\end{array}
\tag{3.69}  \label{3.69}
\end{equation}%
Now, from Lemma 3.4, then it easy to shows that%
\begin{equation}
\begin{array}{l}
\frac{1}{2}\mathbb{E}\left[ h_{xx}\left( x^{\ast }(T),\mathbb{E}(x^{\ast
}(T))\right) x_{1}^{\varepsilon }(T)^{2}\right] \smallskip \smallskip \\ 
\text{ \ \ \ }=\mathbb{E}\dint_{s}^{T}\left\{ \frac{1}{2}H_{xx}(t)\left(
x_{1}^{\varepsilon }(t)\right) ^{2}\right. -\frac{1}{2}Q^{\ast }(t)\sigma
_{y}^{2}(t)\left( \mathbb{E}\left( x_{1}^{\varepsilon }(t)\right) \right)
^{2}\smallskip \smallskip \\ 
\text{ \ \ \ }-\frac{1}{2}Q^{\ast }(t)\left( \delta \sigma (t)\right) ^{2}-%
\frac{1}{2}\dint_{\Theta }Q^{\ast }(t)\left( \delta g(t,\theta \right)
)^{2}\mu (d\theta )\mathbf{I}_{\mathcal{E}_{\varepsilon }}(t)\smallskip
\smallskip \\ 
\text{ \ \ \ }-\frac{1}{2}\dint_{\Theta }\Gamma _{t}^{\ast }(\theta )\left(
\delta g(t,\theta \right) )^{2}\mathbf{I}_{\mathcal{E}_{\varepsilon }}(t)\mu
(d\theta )\smallskip \smallskip \\ 
\text{ \ \ \ }-\left( \mathbb{E}\left( x_{1}^{\varepsilon }(t)\right)
\right) x_{1}^{\varepsilon }(t)\left[ Q^{\ast }(t)f_{y}(t)+Q^{\ast
}(t)\sigma _{x}(t)\sigma _{y}(t)+R^{\ast }(t)\sigma _{y}(t)\right]
\smallskip \smallskip \\ 
\text{ \ \ \ }-\frac{1}{2}\left. \dint_{\Theta }Q^{\ast }(t)\left(
g_{x}\left( t,\theta \right) \right) ^{2}\left( x_{1}^{\varepsilon
}(t)\right) ^{2}\mu (d\theta )\right\} dt+\tau \left( \varepsilon \right) ,%
\end{array}
\tag{3.70}  \label{3.70}
\end{equation}%
using Lemma 3.5\textbf{\ }together with\textbf{\ }(\ref{3.69}) and (\ref%
{3.70}) we obtain

\begin{eqnarray*}
0 &\leq &J^{^{s,\zeta }}\left( u^{\varepsilon }(\cdot )\right) -J^{^{s,\zeta
}}\left( u^{\ast }(\cdot )\right) =-\mathbb{E}\int_{s}^{T}\delta H(t)\mathbf{%
I}_{\mathcal{E}_{\varepsilon }}(t)dt \\
&&-\frac{1}{2}\mathbb{E}\int_{s}^{T}Q^{\ast }(t)\left( \delta \sigma
(t)\right) ^{2}\mathbf{I}_{\mathcal{E}_{\varepsilon }}(t)dt \\
&&-\frac{1}{2}\mathbb{E}\int_{s}^{T}\int_{\Theta }Q^{\ast }(t)(\delta
g\left( t,\theta \right) )^{2}\mathbf{I}_{\mathcal{E}_{\varepsilon }}(t)\mu
(d\theta )dt \\
&&-\frac{1}{2}\mathbb{E}\int_{s}^{T}\int_{\Theta }\Gamma _{t}^{\ast }(\theta
)(\delta g\left( t,\theta \right) )^{2}\mathbf{I}_{\mathcal{E}_{\varepsilon
}}(t)\mu (d\theta )dt+\tau \left( \varepsilon \right) ,
\end{eqnarray*}%
then we get

\begin{eqnarray*}
0 &\leq &J^{^{s,\zeta }}\left( u^{\varepsilon }(\cdot )\right) -J^{^{s,\zeta
}}\left( u^{\ast }(\cdot )\right) =-\mathbb{E}\int_{s}^{T}\delta H(t)dt \\
&&-\frac{1}{2}\mathbb{E}\int_{s}^{T}Q^{\ast }(t)\left( \delta \sigma
(t)\right) ^{2}\mathbf{I}_{\mathcal{E}_{\varepsilon }}(t)dt \\
&&-\frac{1}{2}\mathbb{E}\int_{s}^{T}\int_{\Theta }\left( Q^{\ast }(t)+\Gamma
_{t}^{\ast }(\theta )\right) (\delta g\left( t,\theta \right) )^{2}\mathbf{I}%
_{\mathcal{E}_{\varepsilon }}(t)\mu (d\theta )dt+\tau \left( \varepsilon
\right) .
\end{eqnarray*}%
Finally by using (\ref{2.4}) we deduce%
\begin{eqnarray*}
0 &\leq &\mathbb{E}\int_{s}^{T}\left\{ -H\left( t,x^{\ast },\mathbb{E}\left(
x^{\ast }\right) ,u,\Psi ^{\ast }(t),K^{\ast }(t),\gamma _{t}^{\ast }(\theta
)\right) \right. \\
&&+H\left( t,x^{\ast },\mathbb{E}\left( x^{\ast }\right) ,u^{\ast }(t),\Psi
^{\ast }(t),K^{\ast }(t),\gamma _{t}^{\ast }(\theta )\right) \\
&&-\frac{1}{2}Q^{\ast }(t)\left( \sigma (t,x^{\ast }(t),\mathbb{E}(x^{\ast
}(t)),u)-\sigma (t,x^{\ast }(t),\mathbb{E}(x^{\ast }(t)),u^{\ast
}(t))\right) ^{2}\mathbf{I}_{\mathcal{E}_{\varepsilon }}(t) \\
&&\left. -\frac{1}{2}\int_{\Theta }\left( Q^{\ast }(t)+\Gamma _{t}^{\ast
}(\theta )\right) (g\left( t,x^{\ast }(t),u,\theta \right) -g\left(
t,x^{\ast }(t),u^{\ast }(t),\theta \right) )^{2}\mathbf{I}_{\mathcal{E}%
_{\varepsilon }}(t)\mu (d\theta )\right\} dt \\
&&+\tau \left( \varepsilon \right) .
\end{eqnarray*}%
This completes the proof of Theorem 3.1.$\smallskip \smallskip $

\noindent \textbf{Conclusions. }In this paper, stochastic maximum principle
for optimal stochastic control for systems governed by SDE of mean-field
type with jump processes is proved. The control variable is allowed to enter
both diffusion and jump coefficients and also the diffusion coefficients
depend on the state of the solution process as well as of its expected
value. Moreover, the cost functional is also of Mean-field type. When the
coefficients $f$ and $\sigma $ of the underlying diffusion process and the
cost functional do not explicitly depend on the expected value, \textit{%
Theorem 3.1} reduces to stochastic maximum principle of optimality, proved
in Tang et al., (\cite{tang}, \textit{Theorem 2.1}).$\smallskip \smallskip $

\noindent \textbf{Appendix}\medskip

\noindent The following result gives special case of the It\^{o} formula for
jump diffusions.$\smallskip \smallskip $

\noindent \textbf{Lemma A1. (}\textit{Integration by parts formula for jumps
processes}) Suppose that the processes $x_{1}(t)$ and $x_{2}(t)$ are given
by: for $i=1,2,$ $t\in \left[ s,T\right] :$

\begin{equation*}
\left\{ 
\begin{array}{l}
dx_{i}(t)=f\left( t,x_{i}(t),u(t)\right) dt+\sigma \left(
t,x_{i}(t),u(t)\right) dW(t)\smallskip \smallskip \\ 
\text{ \ \ \ \ \ \ \ \ }+\dint_{\Theta }g\left( t,x_{i}(t_{-}),u(t),\theta
\right) N\left( d\theta ,dt\right) ,\smallskip \smallskip \\ 
x_{i}(s)=0.%
\end{array}%
\right.
\end{equation*}%
Then we get 
\begin{eqnarray*}
\mathbb{E}\left( x_{1}(T)x_{2}(T)\right) &=&\mathbb{E}\left[
\int_{s}^{T}x_{1}(t)dx_{2}(t)+\int_{s}^{T}x_{2}(t)dx_{1}(t)\right] \\
&&+\mathbb{E}\int_{s}^{T}\sigma ^{\ast }\left( t,x_{1}(t),u(t)\right) \sigma
\left( t,x_{2}(t),u(t)\right) dt \\
&&+\mathbb{E}\int_{s}^{T}\int_{\Theta }g^{\ast }\left(
t,x_{1}(t),u(t),\theta \right) g\left( t,x_{2}(t),u(t),\theta \right) \mu
(d\theta )dt.
\end{eqnarray*}

\noindent See Framstad et al., (\cite{framstad}, \textit{Lemma 2.1}) for the
detailed proof of the above Lemma.

\noindent \textbf{Proposition A1.} Let $\mathcal{G}$ be the predictable $%
\sigma -$field on\ $\Omega \times \left[ s,T\right] $, and $f$ be a $%
\mathcal{G}\times \mathcal{B}(\Theta )-$measurable function such that%
\begin{equation*}
\mathbb{E}\int_{s}^{T}\int_{\Theta }\left\vert f\left( r,\theta \right)
\right\vert ^{2}\mu (d\theta )dr<\infty ,
\end{equation*}%
then for all $p\geq 2$ there exists a positive constant $C=C(T,p,\mu (\Theta
))$ such that%
\begin{equation*}
\mathbb{E}\left[ \sup_{0\leq t\leq T}\left\vert \int_{s}^{t}\int_{\Theta
}f\left( r,\theta \right) N(d\theta ,dr)\right\vert ^{p}\right] \leq C%
\mathbb{E}\left[ \int_{s}^{T}\int_{\Theta }\left\vert f\left( r,\theta
\right) \right\vert ^{p}\mu (d\theta )dr\right] .
\end{equation*}%
\noindent \textbf{Proof.} See Bouchard et al., (\cite{bouchard}, \textit{%
Appendix}).$\smallskip \smallskip $

\noindent \textbf{Lemma A2 }(\textit{Martingale representation theorem for
jump processes}). Let $\mathcal{G}$ be a finite-dimensional space and let $%
m(t)$ be an $\mathcal{G-}$valued $\mathcal{F-}$adapted square-integrable
Martingale. Then there exist $q(\cdot )\in \mathbb{L}_{\mathcal{F}%
}^{2}\left( \left[ s,T\right] ,\mathcal{G}\right) $ and $g(\cdot ,\cdot )\in 
\mathbb{M}_{\mathcal{F}}^{2}\left( \left[ s,T\right] ,\mathcal{G}\right) $
such that%
\begin{equation*}
m(t)=m(s)+\int_{s}^{t}q(r)dW(r)+\int_{s}^{t}\int_{\Theta }g(r,\theta
)N(d\theta ,dr).
\end{equation*}

\noindent \textbf{Proof. }See Tang et al., (\cite{tang} \textit{Lemma 2.3}).

\end{document}